\documentclass[11pt]{article}
\usepackage{german}
\usepackage[T1]{fontenc}
\usepackage[ansinew]{inputenc}
\usepackage{latexsym}
\usepackage{textcomp}
\usepackage{color}
\usepackage{amscd}
\usepackage{amsmath}

\input{amssym.def}
\input{amssym}

     \newcommand{\fa}{\goth{a}}
     \newcommand{\fb}{\goth{b}}
     \newcommand{\fc}{\goth{c}}

     \newcommand{\fu}{\goth{u}}

     \newcommand{\fF}{\goth{F}}
     \newcommand{\fG}{\goth{G}}
     \newcommand{\fH}{\goth{H}}

     \newcommand{\fO}{\goth{O}}

     \newcommand{\fV}{\goth{V}}

     \newcommand{\F}{\Bbb{F}}
     
     \newcommand{\Ll}{\Bbb{L}}  

     \newcommand{\Q}{\Bbb{Q}}
     
     \newcommand{\Z}{\Bbb{Z}}

    \newcommand{\qed}{{\hfill$\Box$}}

    \newcommand{\ol}[1]{\overline{#1}}

    \newcommand{\ti}[1]{\tilde{#1}}
    \newcommand{\group}[1]{\langle{#1}\rangle}

    \newcommand{\ot}{\otimes}
    \newcommand{\me}{^{-1}}
    \newcommand{\mal}{^{\times}}
    
    \newcommand{\df}{\stackrel{\mathrm{def}}{=}}
    \newcommand{\mr}{\mathrm}
    \newcommand{\clo}{^{\mr{c}}}

    \newcommand{\qg}{{\Bbb{Q}G}}

    \newcommand{\zl}{{\Bbb{Z}_l}}
    
    \newcommand{\ql}{{\Bbb{Q}_l}}

    \newcommand{\into}{\rightarrowtail}
    \newcommand{\onto}{\twoheadrightarrow}
    \newcommand{\lto}{\longrightarrow}
    \newcommand{\da}{\downarrow}

    \newcommand{\pl}{\parallel}
    \newcommand{\pht}{\phantom}

    \def\daz#1{#1\da\pht{#1}}

    \newcommand{\sda}{\da_{\raisebox{-1.4mm}{$\hsp{-1}\check{}$}}}

    \newcommand{\ka}{\kappa}
    \newcommand{\ga}{\gamma}
    \newcommand{\Ga}{\Gamma}
    \newcommand{\si}{\sigma}
    \newcommand{\Si}{\Sigma}
    
    \newcommand{\De}{\Delta}
    
    \newcommand{\la}{\lambda}
    \newcommand{\La}{\Lambda}
    \newcommand{\be}{\beta}
    
    \newcommand{\om}{\omega}
    \newcommand{\al}{\alpha}

    \newcommand{\noi}{\par\noindent}
    \newcommand{\sn}{\par\smallskip\noindent}
    \newcommand{\mn}{\par\medskip\noindent}
    \newcommand{\bn}{\par\bigskip\noindent}
    \newcommand{\bbn}{\par\bigskip\bigskip\noindent}
    \newcommand{\bbbn}{\par\bigskip\bigskip\bigskip\noindent}
    \newcommand{\Section}[2]{\bbn {\large #1\,. \ {\sc #2}}
                             \nopagebreak
                             \nz}

    \newcommand{\nf}[2]{\\[1.5ex]
                        \bmp{1cm}
                         (#1)
                        \emp 
                        \bmp{13.5cm}
                         \bct
                          $#2$
                         \ect
                        \emp\\[1.5ex]
         }

    \newcommand{\nz}{\\[1ex]}

    \newcommand{\hsp}[1]{\hspace*{#1mm}}

    \newcommand{\mmargin}{
     \textheight 230truemm
     \textwidth 155truemm
     \topmargin -10truemm
     \oddsidemargin 5truemm
     \evensidemargin 5truemm
     }

    \newcommand{\bmp}{\begin{minipage}}
    \newcommand{\emp}{\end{minipage}}
    \newcommand{\btb}{\begin{tabular}}
    \newcommand{\etb}{\end{tabular}}
    \newcommand{\barr}{\begin{array}}
    \newcommand{\earr}{\end{array}}
    \newcommand{\bit}{\begin{itemize}}
    \newcommand{\eit}{\end{itemize}}
    \newcommand{\ben}{\begin{enumerate}}
    \newcommand{\een}{\end{enumerate}}
    \newcommand{\bct}{\begin{center}}
    \newcommand{\ect}{\end{center}}
    \newcommand{\bfr}{\begin{flushright}}
    \newcommand{\efr}{\end{flushright}}
    \newcommand{\bea}{\begin{eqnarray*}}
    \newcommand{\eea}{\end{eqnarray*}}
    \newcommand{\bqo}{\begin{quote}}
    \newcommand{\eqo}{\end{quote}}
    \newcommand{\bdc}{\begin{description}}
    \newcommand{\edc}{\end{description}}
    \newcommand{\bdia}{\begin{CD}}
    \newcommand{\edia}{\end{CD}}

    \definecolor{light}{gray}{.3}

    \newcommand{\Hom}{\mathrm{Hom}}

    \newcommand{\tr}{\mathrm{tr\,}}

    \newcommand{\ind}{\mathrm{ind\,}}
    \newcommand{\infl}{\mathrm{infl\,}}

    \newcommand{\res}{\mathrm{res\,}}
    
    \newcommand{\Det}{\mathrm{Det\,}}
    \newcommand{\sr}[2]{{\,\stackrel{#1}{#2}\,}}
    \newcommand{\fra}[2]{{\,\frac{#1}{#2}\,}}

    \newcommand{\theorem}{\sn
           \bdc
           \item[{\sc Theorem.}] \em }
    \newcommand{\Theorem}[1]{\sn
           \bdc
           \item[{\sc Theorem {#1}.}]  \em }

    \newcommand{\Stop}{\edc \sn\rm} 

    \newcommand{\Lemma}[1]{\sn
           \bdc
           \item[{\sc Lemma {#1}.}] \em }
    
    \newcommand{\Proposition}[1]{\sn
           \bdc
           \item[{\sc Proposition {#1}.}] \em }

    \newcommand{\proof}{{\sc Proof.} \ }
    
    \newcommand{\Remark}[1]{\sn{\sc Remark {#1}. \ }}

\mmargin

\def\lw{{\La_\wedge}}
\def\lwa{{\lw A}}
\def\lba{{\La_\bullet A}}

\def\lg{{\La G}}
\def\qq{{\mathcal{Q}}}
\def\qg{{\qq G}}
\def\qwg{{\qq_\wedge G}}
\def\lwg{{\lw G}}
\def\rlg{{R_lG}}

\def\HOM{{\mr{HOM}}}
\def\Det{{\mr{Det}}}
\def\defl{{\mr{defl}}}
\newcommand{\LL}{{{\mbox{{\boldmath$\mr{L}$}}}}}
\def\Tr{{\mr{Tr}}}
\def\tr{{\mr{tr}}}
\def\infl{{\mr{infl}}}
\def\ver{{\mr{ver}}}
\def\qwcgak{{\qq_\wedge\clo\Ga_k}}
\def\ab{{\mr{ab}}}
\def\Res{{\mr{Res}}}
\def\gab{{G^\ab}}
\def\uab{{U^\ab}}
\def\ggc{{[G,G]}}
\def\uuc{{[U,U]}}
\def\gq{{\ol G}}
\def\aq{{\ol A}}
\def\nq{{\ol N}}
\def\uq{{\ol U}}
\def\gs{{G'}}

\def\supp{{\mr{supp}}}
\def\st{{\mr{stab}}}

\def\ep{{\epsilon}}

 \def\muq{{\mu_Q}}
 \def\vt{{\vartheta}}
\def\ue{{U_1}}

\def\sign{{\mr{sgn}}}
\def\vr{{\mr{vr}}}

\def\muq{{\mu_Q}}
\def\muqs{{\mu_{Q'}}}
\def\stab{{\mr{stab}}}
\def\symq{{\mr{Sym}(Q)}}
\def\symqv{{\mr{Sym}(Q/V)}}
\def\kas{{\ka(s)}}
\def\kase{{\ka(s_1)}}
\def\kasz{{\ka(s_2)}}
\def\kass{{\ka(\si s)}}
\def\kasme{{\ka(s)\me}}
\def\sqv{{s\in Q/V}}
\def\vv{{v\in V}}
\def\ska{{\si_\ka}}
\def\qw{{\qq_\wedge}}
\def\lwb{{\La^\be_\wedge}}
\def\qwb{{\qq^\be_\wedge}}
\def\lwc{{\La_\wedge\clo}}
\def\qwc{{\qq_\wedge\clo}}
\def\dal{{\sr{|}{\da}}}

\begin{document}
 \title{On the `main conjecture' of equivariant Iwasawa theory}
 \author{Jürgen Ritter \ $\cdot$ \ Alfred Weiss \
 \thanks{We acknowledge financial support provided by DFG and NSERC.}
 }
\date{\pht\today}

 \maketitle
 \noi The `main conjecture' of equivariant Iwasawa theory concerns the
 situation where

 \ben\item[]  $l$ is a fixed odd prime number\,,
  $K/k$ is a Galois extension of totally real
 number fields with $k/\Q$ and $K/k_\infty$ finite, where $k_\infty/k$ is
 the cyclotomic $\zl$-extension (we set $G=G(K/k)$ and $\Ga_k=G(k_\infty/k)$)\,,
 \item[] $S$ a fixed finite set (which will normally be suppressed in the notation) of primes of $k$ containing
 all primes which ramify in $K$ and all archimedean primes\,,
 and $M$ is the maximal abelian $l$-extension of $K$ unramified outside $S$ (we set $X=G(M/K)$)\,.
 \een
  It asserts that a canonical refinement $\mho=\mho_S$ of the Iwasawa module $X$ is determined by the Iwasawa $L$-function $L_{K/k}=L_{K/k,S}$.
 The data $\mho$ and $L_{K/k}$ have been defined in [RW2] and we briefly recapitulate what is needed here.

 \mn Denote by $\lg\,,\,\qg$ the completed
 group ring $\zl[[G]]$ of $G$ over $\zl$ and its total ring of fractions, respectively.
 The localization sequence of $K$-theory
 $$K_1(\La G)\to K_1(\qg)\sr{\partial}{\to} K_0T(\La G)\to K_0(\qg)$$
 has $\mho$ in $K_0T(\La G)$ (see [loc.cit., p.563]).

 \mn The reduced norms of the Wedderburn components of the
 semi-simple algebra $\qg$ induce the map (see [loc.cit., Theorem 8])
 $$\Det:K_1(\qg)\to \Hom^\ast(\rlg,(\qq\clo\Ga_k)\mal)\,.$$
 The Iwasawa $L$-function $L_{K/k}$ is derived from the $S$-truncated $l$-adic Artin $L$-functions (see [loc.cit., Proposition 11]) and
 belongs to the above group $\Hom^\ast$.

 \mn The equivariant `main conjecture' asserts that \ {\em there is a unique element $\Theta\in K_1(\qg)$
 satisfying $\Det(\Theta)=L_{K/k}$ and, moreover, that this $\Theta$ has
 $\partial(\Theta)=\mho$\,.} In other words, $\Theta$ is the nonabelian pseudomeasure of $K/k$.
\nopagebreak[4]

 \theorem If Iwasawa's $\mu$-invariant $\mu_{K/k}$ vanishes, then the `main conjecture' of equivariant Iwasawa theory for $K/k$ holds,
 up to its uniqueness statement.\Stop
 There are now much more general conjectures on nonabelian Iwasawa theory, especially in [FK] and [Ka], which also have extensive
 bibliographies.

 \mn The proof of the {\sc Theorem} depends on our previous work. So we start by describing certain reductions in terms of
 Propositions for which we provide programmatic proofs in \S1.

 \mn Our first reduction, in [RW3,4], of the proof of the {\sc Theorem} loses the uniqueness assertion of the conjecture; we accepted this
 because uniqueness would follow from $SK_1(\qg)=1$, as conjectured by Suslin (see [RW2, Remark E]). It also requires us
 to assume that $\mu_{K/k}=0$\,, as conjectured by Iwasawa,
 in order to bring in the localization $\La_\bullet G$ of $\La G$ obtained by inverting
 all central elements which are regular modulo $l$ and even the completion $\lwg$ of the localization, which permits a logarithmic approach yielding

 \Proposition{1} The {\sc Theorem} holds if, and only if, $L_{K/k}\in \Det K_1(\lwg)$\,.\Stop
 Note that $\mho$ no longer appears: this is a consequence of the Main Conjecture of classical Iwasawa theory, proved by Wiles [Wi].

 \mn Define $K/k$ to be $l$-{\em elementary}\,, if
 $G(K/k)=\group{z}\times G[l]$ is a direct product of a finite
 cyclic group $\group{z}$ of order $|z|$ prime to $l$ and a
 pro-$l$ group $G[l]$.

 \mn The second reduction is due to induction techniques, [RW4],
 for $l$-adic characters with open kernel, from which we obtain

 \Proposition{2} If $L_{K/k}\in\Det K_1(\lw G(K/k))$ holds
 whenever $K/k$ is $l$-elementary, then $L_{K/k}\in\Det K_1(\lw
 G(K/k))$ for arbitrary extensions $K/k$. \footnote{Actually, the uniqueness assertion
 of the `main conjecture' can also be reduced to $l$-elementary
 $K/k$, see [La, Chapter 4].}\Stop
 {\em From now on, $K/k$ is always $l$-elementary} \footnote{if
 not otherwise implied}.

 \mn The third reduction comes from a logarithmic interpretation of `\,$L_{K/k}\in\Det K_1(\lwg)$\,' by means of diagram (D1) in \S1\,:
 the point is that the nontriviality of the kernel of `Det' on $K_1(\lwg)$ is an obstacle to inductive constructions of preimages of the
 Iwasawa $L$-function.

 \Proposition{3}\quad $L_{K/k}\in \Det K_1(\lwg)\iff t_{K/k}\in T(\lwg)$\Stop
 Here, the {\em logarithmic pseudomeasure} $t_{K/k}\in T(\qwg)$ is recalled in \S1. The proposition is a consequence of the validity of the
 {\em torsion congruence}\,, which follows from the $q$-expansion principle of Deligne and Ribet
 [DR]. The proof, [RW5,\S3], of the proposition for pro-$l$ extensions is
 extended to $l$-elementary extensions in \S1.

 \mn The methods developed for proving Proposition 3 combined with its emphasis on working in the $T$-world motivates a new ingredient,
 the restriction $\Res_G^U:T(\qwg)\to T(\qw U)$ for open subgroups $U\ge\group{z}$ of $G$, and to

 \Proposition{4} If $G$ has an abelian subgroup $A$ of index $l$, then the {\sc Theorem} holds \footnote{The generalization of
 [RW8, Theorem] to the $l$-elementary case is already in [La, Chapter 2].}.\Stop
 In fact, most of the ideas of the present paper have already appeared in the proof of Proposition 4, in embryo. Thus the proof of
 the {\sc Theorem} can be viewed as
 a full generalization of [RW8], i.e., we {\em generalize} the Wall
 congruence
 (see Proposition 5), the torsion congruence (see Proposition 6) and the proof of Proposition 4 (see Theorem 7).

 \mn Let $\mu_Q$ denote the Möbius function of the partially ordered set of subgroups of the finite $l$-group
 $Q$. Recall that $\mu=\muq$ is defined by $$\mu(1)=1\,,\ \mu(Q')=-\sum_{1\le Q''<Q'}\mu(Q'')\quad\mr{for}\quad1\neq Q'\le Q\,.$$

 \Proposition{5} {\rm (``Möbius-Wall'')} \quad Let $A$ be an abelian normal open subgroup of $G=G(K/k)$ so that $Q=G/A$ is a
 finite $l$-group. If $\ep$ is a unit of $\lwg$,
 then $$\sum\limits_{A\le U\le G}\mu_Q(U/A)\ver_U^A(\res_G^U\ep)\equiv0\mod\tr_Q(\lwa)\,.$$ \Stop
 Here, $\ver_U^A:\lw U\to\lwa$ extends the group transfer $U\to \uab\to A$ to a ring homomorphism between their Iwasawa algebras in the customary way
  \footnote{The inverse system of $\zl$-linear extensions $\zl[U/V]\to\zl[A/V]$ of
  the group tranfers $U/V\to A/V$, where $V$ runs through the normal open subgroups of $U$ contained in $A$ ,
   give rise to the transfer $\La U\to\La A$, which is a ring homomorphism
   that can be localized and completed. Note that $\ver_U^A=\ver_\uab^A\defl_U^\uab$ factors through $\uab$.}.
 Proposition 5 is proved in \S2.

 \mn Let $S'$ denote the set of all non-archimedean primes of $S$ and $k\subseteq f\subset K$ with $[f:k]<\infty$. The pseudomeasure $\la_f=\la_{f,S'}$
 of [Se] is associated to the maximal abelian $S'$-ramified extension $f_{S'}$ of $f$. We extend this notation to intermediate fields
 $F$ of $f_{S'}/f_\infty$ by defining  $\la_{F/f}=\defl_{G(f_{S'}/f)}^{G(F/f)}\la_f$\,, where $\defl_{G(f_{S'}/f)}^{G(F/f)}$ is the
 deflation map $\lw G(f_{S'}/f)\to\lw G(F/f)$. Now we can state

 \Proposition{6} Notation as in Proposition 5, $$\sum\limits_{A\le U\le G}\mu_Q(U/A)\ver_\uab^A(\la_{K^{[U,U]}/K^U})\equiv0\mod\tr_Q(\lba)\,.$$\Stop
 The proof of the Proposition follows from [RW9, Theorem], which again builds on [DR]; see \S3.

 \mn As has been mentioned in the introduction of [RW9], combining the `main conjecture' with Proposition 5 implies Proposition 6. Conversely,
 fusing Propositions 5 and 6 leads to a partial generalization of Proposition 4 and its proof.
 \Theorem{7} Let $A\ge\group{z}$ be an abelian normal open subgroup of $G=G(K/k)$ and $C$ a central subgroup
 of exponent $l$ contained in $A$. If $t_{K^C/k}$ is integral, then there exists a $\xi\in T(\lwg)$
 with $\defl_G^{G/C}\xi =t_{K^C/k}$ and $\Res_G^A\xi =t_{K/K^A}$\,.\Stop
 Parts of the proof of Theorem 7, in \S4, are easier versions of arguments appearing already in [RW8] but
 repeated here for the convenience of the reader. The real strength of the theorem is that it serves as
 a catalyst for the proof of the {\sc Theorem}.

 \mn The {\sc Theorem} is proved in \S5 by making suitable modifications in $T(\lwg)$ of the element $\xi$ provided by Theorem 7.

 \mn Section 6 contains the necessary extension of the integral
 logarithm to $l$-elementary groups. This is based on using
 projections to the integral logarithm for pro-$l$ groups with
 unramified coefficients, which is already in [RW3]. It also
 discusses `Res' and the $l$-elementary ingredients of the proof of Proposition 3.

 \bn Finally, in a short appendix, we take the opportunity of
 straightening out an inaccuracy in the proof of [RW2,
 Proposition 12]. In it we have referred to [RWt] where, however,
 Leopoldt's conjecture is assumed to hold. In the appendix we now
 outline an argument which is not based on this conjecture.

 \Section{1}{``Proofs'' of Propositions 1,2,3,4}
  Proposition 1, for pro-$l$ groups, is reduced via the Main Conjecture, as in [RW3,\S1], to the stronger form of Theorem B in [loc.cit.,\S6].
  The essential ingredient
  in this form of Theorem B is the integral group logarithm $\Ll$ defined by the commutative square (see [loc.cit., Proposition 11])
 \nf{D1}{\barr{cccc}
  K_1(\La_\wedge G)&\sr{\Ll}{\to}&T(\qq_\wedge G)&\\
 \daz{\Det}&&\daz{\Tr}\hsp{-4}\simeq&\\
 \mr{HOM}(\rlg,(\La_\wedge\clo\Ga_k)\mal)&\sr{\LL}{\to}&\Hom^\ast(\rlg,\qwcgak)&,\earr}
 with $T(R)=R/[R,R]$ for any ring $R$, where $[R,R]$ is the additive subgroup generated by all Lie commutators $[a,b]=ab-ba\,,\,a,b\in R$,
 and with the isomorphism `Tr' induced by the reduced trace of $\qwg$. For HOM and $\LL$ see [loc.cit.,p.37]. The logarithm $\Ll$ is called
 integral because it takes values in $T(\lwg)$. Note that the Wall congruence, which plays an important role later, makes its
 first appearance in [loc.cit., Lemma 12].
 It is also important to observe, [loc.cit.,p.42], that \footnote{$L_{K/k}(\chi)^l\equiv\Psi(L_{K/k}(\psi_l\chi))\mod l\lw(\Ga_k)$}  $L_{K/k}\in\HOM$.

 \sn The generalization to arbitrary extensions $K/k$ is carried out in [RW4, Theorems (A) and (B)].\qed

\mn Proposition 2 is [RW4, Theorem (C)]. \qed

 \mn We should stress that (D1) is available to define $\Ll$ for
 arbitrary groups $G=G(K/k)$ and that its integrality property,
 $\Ll(K_1(\lwg))\subseteq T(\lwg)$, holds when $G$ is
 $l$-elementary. This will be discussed in \S6.

 \mn As a direct consequence, there is a unique element $t_{K/k}\in
T(\qwg)$, the
 logarithmic pseudomeasure of $K/k$, such that $\Tr(t_{K/k})=\LL(L_{K/k})$\,. By (D1) and i) of Lemma {\sc i}
 in \S6, $L_{K/k}\in\Det K_1(\lwg)$ implies $t_{K/k}\in T(\lwg)$,
 which is the easy implication in Proposition 3.

 \mn Proposition 3 for pro-$l$ groups is proved in [RW5] and [RW7]. The converse direction,
 [RW5, Proposition 2.4] yields
 $w L_{K/k}\in\Det K_1(\lwg)$ where
 $w$ is the unique
 torsion element in $\HOM(\rlg,(\La_\wedge\clo\Ga_k)\mal)$ deflating to 1 on applying  $\defl_G^\gab$; this is extended to the $l$-elementary case in i) of Lemma {\sc k}.
 Then [RW5, Theorem, p.1096] reduces the question of whether
 $L_{K/k}\in\Det K_1(\lwg)$ to the analogous one for Galois extensions $K/k$ with $G=G(K/k)$ having an abelian subgroup of index $l$\,: this carries
 over to the $l$-elementary case without change. Hence in
 order to deduce $w=1$ we may thus assume that $G$ has an abelian
 subgroup $G'$ of index $l$, and then apply the extension in Lemma {\sc k}\,,\,ii) of [loc.cit., Proposition
 3.2] to $l$-elementary groups, to obtain the equivalence
$$w=1\iff\ver(\la_{K^\ggc/k})\equiv\la_{K/K^{G'}}\mod\tr_{G/G'}(\La_\wedge
 G')\,,$$ in which the displayed congruence is referred
 to earlier as the torsion congruence. Here, $\la_{F/f}$ is the pseudomeasure for the extension $F/f$ and
 $\ver:\lw\gab\to\lw G'$ is induced by the group
 transfer $\gab\to G'$\,. Note that the proof of Proposition 3.2 depends on the Wall congruence mentioned above.
 This torsion congruence is proved in [RW7] (or [RW9]) by interpreting the methods of [DR] on the Galois side as in [Se].
 We point out right away that the validity of the torsion
 congruence persists in Proposition 6 because every open subgroup of index $l$ in $G$
 contains $\group{z}$; since $K\subset L_{S'}^+$ (in the notation of
 the first paragraph of \S3) we can specialize the group
 $H_{S'}^+$ to our $A$.\qed

 \Remark{A} Consider the localization $\zl[[H_{S'}^+]]_\bullet$ of [RW7, bottom of p.715], which results from inverting the multiplicative set of elements of $\zl[[H_{S'}^+]]$ whose image in $\zl[[\Ga_L]]$ is not in $l\zl[[\Ga_L]].$ Observe that if $A\ge\group{z}$ and $z\neq1$, then $\hat z=\sum_{i=0}^{|z|-1}z^i$ has image $|z|$ in $\zl[[\Ga_L]]$ and so is not in $l\zl[[\Ga_L]]$\,, hence $\hat z\in\zl[[H_{S'}^+]]_\bullet\mal$ and $\hat z(z-1)=0$  imply $z=1$ in $\zl[[H_{S'}]]_\bullet$. It follows that specializing $H_{S'}^+$ to $A$ induces a map from $\zl[[H_{S'}^+]]_\bullet$ to the current $\lba$ only when $A$ is a pro-$l$ group.

  \Remark{B} It should perhaps be added that  $\La_\bullet G=\Si\me\lg$ with $\Si=\La\Ga\setminus l\cdot\La\Ga$ for any central open subgroup $\Ga\simeq\zl$ of an arbitrary $G=G(K/k)$. Note that $\Si\me\La\Ga$ has the unique maximal ideal $l\Si\me\La\Ga$.
 So it suffices to show that every element $c\in\Si\me\La G$, which is (left) regular modulo $l$, is a unit of $\Si\me\La G$.
 For this consider right multiplication $\Si\me\La G\sr{\cdot c}{\lto}\Si\me\La G$ by $c$.   Since $\Si\me\La G/l$ is a finite
 dimensional $\Si\me\La\Ga/l$-vector space, $c\mod l$ has a (left) inverse in $\Si\me\La G/l $, hence $c$ a left inverse  $b$ in $\Si\me\La G$
 by Nakayama's lemma. Since $b\mod l$ is now also (left) regular modulo $l$, the same argument provides $a\in\Si\me\La G$ with
  $ab=1$. Then $a=abc=c$, so $c$ is a unit.

 \bn Proposition 4, for pro-$l$ groups, is shown in [RW8]. The restriction $\Res_G^U$ is motivated by the
 left square of the commutative diagram
 \nf{D2}{\barr{cccccc} K_1(\lwg)&\sr{\Ll}{\to}&
 T(\qwg)&\sr{\Tr}{\to}&\Hom^\ast(\rlg,\qwcgak)&\\
 \daz{\res_G^U}&&\daz{\Res_G^U}&&\daz{\Res_G^U}&\\
 K_1(\lw U)&\sr{\Ll}{\to}& T(\qq_\wedge
 U)&\sr{\Tr}{\to}&\Hom^\ast(R_lU,\qq_\wedge\clo\Ga_{K^U})&\earr}
 and discussed in detail in [loc.cit., \S1 and Appendix]. Its extension to $l$-elementary groups with $U\ge\group{z}$ is again in \S6; see Lemma {\sc j}.\qed

 \Section{2}{Proof of Proposition 5}
 Fix a set of coset representatives  $r_q$ of $A$ in $G$,
 whence \
 $G=\dot{\bigcup}_{q\in Q}r_qA\,,\,q=r_qA$ and $r_{q_1}r_{q_2}=r_{q_1q_2}a_{q_1,q_2}$ with $a_{q_1,q_2}\in
 A$ a 2-cocycle, so $a_{q_1,q_2q_3}a_{q_2,q_3}=a_{q_1q_2,q_3}a_{q_1,q_2}^{q_3}\,.$
 Further, let $\Si=\symq$ denote the symmetric group on the elements of
 $Q$. It carries the natural (right) $Q$-action
 $$\pi^q(q_1)=\pi(q_1q\me)q\,,\quad q,q_1\in Q\,,\,\pi\in\Si$$ satisfying $(\pi_1\pi_2)^q=\pi_1^q\pi_2^q$\,.
 For $V\le Q$, the set of fixed points $\Si^V$ of $V$ in
 $\Si$ is thus a subgroup of $\Si$. Note that
 $\pi\in\Si^V$ has $\pi(qv)=\pi(q)v$ for all $q\in Q\,,\,v\in V$.

 \Lemma{a} Let $U$ be a subgroup of $G$ containing $A$, set $V=U/A$ and fix a section $\ka:Q/V\to Q$\,, so $(\ka s)V=sV$ for $s\in Q/V$. Let $\ep=\sum_{q\in Q}r_qe_q$\,, with $e_q\in\lwa$, be a unit in $\lwg=\bigoplus_{q\in Q}r_q\cdot\lwa$\,. Then $$\ver_U^A\res_G^U\ep=\sum_{\pi\in\Si^V}\sign(\pi)\prod_{q\in
 Q}a_{\pi(q)q\me,q}\prod_{s\in Q/V}\ver_U^A(e_{\pi(\ka s)\kasme}^{\kas})\
 .$$\Stop
 \proof Writing $\lwg=\bigoplus_{s\in Q/V}r_\kas\lw U$, then $$\ep r_\kase=\sum_{s_2\in Q/V}r_\kasz\Big(\sum_\vv r_va_{\kasz,v}\me a_{\kasz v\kase\me,\kase}\cdot e_{\kasz v\kase\me}^\kase\Big)$$ (with the term in parenthesis in $\lw U$). The ring homomorphism $\ver_U^A:\lw U\to\lwa$ induces the map $\ver_U^A:K_1(\lw U)\to K_1(\lwa)\sr{\det}{=}(\lwa)\mal$ and we compute $\ver_U^A\res_G^U\ep$ by applying $\ver_U^A$ to the matrix of the action of $\ep$ on the right $\lw U$-module $\lwg$ to get
 $$\barr{l}\ver_U^A\res_G^U\ep= \sum\limits_{\si\in\symqv}\sign(\si)\prod\limits_\sqv\Big(\sum_\vv\ver_U^A(r_va_{\kass,v}\me a_{\kass v\kasme,\kas}\cdot e_{\kass v\kasme}^\kas)\Big)\\
 =\sum\limits_{\si\in\symqv}\sign(\si)\sum\limits_{f:Q/V\to V}\prod\limits_\sqv\ver_U^A\Big(r_{f(s)}a_{\kass,f(s)}\me a_{\kass f(s)\kasme,\kas}\cdot e_{\kass f(s)\kasme}^\kas\Big)\earr$$ where $f$ varies over all functions $Q/V\to V$\,, hence
 \bn(2.1)\hsp{20}$\ver_U^A\res_G^U\ep=$ {\small
 $$\sum\limits_{\si\in\symqv}\sign(\pi)\sum\limits_f\Big(\prod\limits_\sqv\prod\limits_\vv a_{f(s),v}a_{\kass,f(s)}^{-v}a_{\kass f(s)\kasme,\kas}^v\Big)\prod\limits_\sqv\ver_U^A(e_{\kass f(s)\kasme}^\kas)$$}
 \sn because $\ver_U^A(r_{f(s)})=\prod_\vv a_{f(s),v}$ and $\ver_U^Aa=\prod_\vv a^v$ for $a\in A$.

 \mn We next simplify the above double product for a fixed $f:Q/V\to V$ to get
 \nf{2.2}{\prod_\sqv\prod_\vv a_{\kass f(s)\kasme,\kas v}\,.}
 Namely, the cocycle relation for the triple $(\kass, f(s), v)$\,, $$a_{f(s),v}a_{\kass,f(s)}^{-v}=a_{\kass, f(s)v}\me a_{\kass f(s),v}$$ turns the double product into $$\prod\limits_\sqv\prod\limits_\vv a_{\kass,f(s)v}\me a_{\kass f(s),v}a_{\kass f(s)\kasme,\kas}^v\,,$$
 and from the triple $(\kass f(s)\kasme,\kas,v)$ it then becomes $$\prod_\sqv\prod_\vv a_{\kass, f(s)v}\me a_{\kas,v}a_{\kass f(s)\kasme,\kas v}\,.$$
 Now the substitutions\quad $v\rightsquigarrow f(s)\me v\,,\,s\rightsquigarrow\si\me(s)$\quad yield $$\prod_\sqv\prod_\vv a_{\kass,f(s)v}=\prod_\sqv\prod_\vv a_{\kas,v}\quad\mr{confirming}\ (2.2)\,.$$
 We continue by reparametrising the maps $f:Q/V\to V$ in (2.1) in terms of the kernel $\Si_0^V$ of the group homomorphism
 $$\Si^V\sr{\sim}{\to}\symqv\,,\,\pi\mapsto\ti\pi\,,\,\ti\pi(qV)=\pi(q)V\,.$$
 {\sc Claim\,:}{\em \ben\item For every $\si\in\symqv$ there is a unique $\ska\in\Si^V$ with $\ska\ka=\ka\si$. The map $\si\mapsto\ska:\symqv\to\Si^V$ is a  group homomorphism splitting $\sim$.
 \item There is a bijection $\tau\leftrightarrow f$ between $\Si_0^V$ and $\{f:Q/V\to V\}$ given by \ $\tau(\kas v)=\kas f(s)v\,,\,f(s)=\kasme\tau(\kas)$\,.
\item $\sign(\ska)=\sign(\si)^{|V|}=\sign(\si)$ and
$\sign(\tau)=1$ for all $\tau\in \Si_0^V$\,.\een}

\mn\proof If
$\ska\in\Si^V$ exists, then $\ska(q)=\ska(\ka(qV)\ka(qV)\me q)=(\ska\ka)(qV)\ka(qV)\me
q=\ka(\si(qV))\ka(qV)\me q$ as $\ka(qV)\me q\in V$\,; and
conversely. Finally, $\sign(\ska)=\sign(\si)^{|V|}$ since $\{s_i:i\mod b\}$ a cycle of $\si$ implies $\{\ka(s_i)v:i\mod b\}$ is a cycle of $\ska$ for each $v\in V$; and
$\sign(\tau)=$ $\prod_\sqv\sign(\tau\ \mr{on}\ \kas
V)=\prod_\sqv\sign (v\mapsto f(s)v\ \mr{on}\
V)=\prod_\sqv(-1)^{|V|-[V:\group{f(s)}]}=1\,.$ \qed

\mn Note that every $q\in Q$ is a unique product $q=\kas v$ with
$s\in Q/V\,,\,v\in V$\,, and every $\pi\in\Si^V$ is a unique
product $\pi=\ska\tau$ with $\si\in\symqv\,,\,\tau\in\Si_0^V$\,.
Substituting (2.2) for the double product in (2.1) and using
2.,3.~of the {\sc Claim}\,,
$$\barr{l}\kass f(s)\kasme=\ska(\kas) f(s)\kasme=\ska(\kas f(s))\kasme\\
=\ska(\tau(\kas))\kasme=(\ska\tau)(\kas)\kasme=(\ska\tau)(\kas
v)(\kas v)\me=\pi(q)q\me\ ,\earr$$ we obtain the assertion of
Lemma {\sc a}.\qed

 \bbn If, as before, $A\le U\le G$ has $V=U/A$, then, for
 $e\in\lwa$, \ $\vr_V(e)\df\ver_U^A(e)$ \ defines a ring
 endomorphism of $\lwa$ satisfying $\vr_V(a)=\prod_{v\in
 V}a^v$ for all $a\in A$. This condition determines $\vr_V$ uniquely, as $A$ `generates'
 $\lwa$ additively\,: indeed, picking a central open
 $\Ga\simeq\zl$ in $A$ and writing
 $A=\dot{\bigcup}_aa\Ga\,,\,\lwa=\bigoplus_aa\cdot\lw\Ga$\,, the
 element $e$ becomes $e=\sum_aac_a$ for suitable $c_a\in\lw \Ga$
 and $\vr_V(e)=\sum_a\vr_V(a)\Psi(c_a)$ with
 $\Psi:\lw\Ga\to\lw\Ga$ the ring homomorphism induced by
 $\ga\mapsto\ga^{|V|}$ for $\ga\in\Ga$. In particular, we have
 \nf{2.3}{\vr_V(e)^q=\vr_{V^q}(e^q)$ for $V\le
 Q\,,\,e\in\lwa\,,\,q\in Q\,,}
 because
$\vr_V(a)^q=\prod_{v\in V}a^{vq}=\prod_{v\in
 V}a^{qv^q}=\prod_{w\in V^q}a^{qw}=\vr_{V^q}(a^q)$ for all $a\in A$.

\Lemma{b} For all $Q'\le Q$ and all $e\in\lwa$\,, \
$\sum\limits_{V\le Q'}\muqs(V)\prod\limits_{s\in
 Q'/V}\vr_V(e^s)\equiv0\mod\tr_{Q'}(\lwa)$\,, \\ where $s\in Q'/V$ now means $Q'=\dot{\bigcup}_ssV$.\Stop
 \proof
 We first observe that this holds for all $a\in A$
 because $\vr_V(a^s)=\prod_{v\in V}a^{sv}$\,, so $\prod_{s\in
 Q'/V}\vr_V(a^s)=\prod_{q'\in Q'}a^{q'}$ is independent of $V\le Q'$, and
 $\sum_{V\le Q'}\muqs(V)=0$.
 \sn It therefore suffices to prove
 additivity of the left side of the claimed congruence, i.e.,
 $$\barr{l}\sum_{V\le Q'}\muqs(V)\prod_{s\in
 Q'/V}\vr_V((e_0+e_1)^s)\\ \equiv\sum_{V\le Q'}\muqs(V)\prod_{s\in
 Q'/V}\vr_V(e_0^s)+\sum_{V\le Q'}\muqs(V)\prod_{s\in
 Q'/V}\vr_V(e_1^s)\mod\tr_{Q'}(\lwa)\ .\earr$$

 \mn We proceed by induction on $|Q'|$; the case $Q'=1$ is trivial.

 \sn Let $F=F(Q')$
 denote the set of maps $f$ from $Q'$ to $\F_2$, with $Q'$-action
 $(fq')(x)=f(x(q')\me)$ for all $x\in Q'$. Then
 $$\prod_{s\in Q'/V}(e_0+e_1)^s=\sum_{f\in F^V}\prod_{s\in
 Q'/V}e_{f(s)}^s\,,$$ because the set of fixed points $F^V$ of $V$ on $F$
 is the set of all $f:Q'/V\to\F_2$\,. Defining $\fF=\fF(Q')=\{(V,f):V\le
 Q'\,,\,f\in F^V\}$ we have
 $$\sum_{V\le Q'}\muqs(V)\prod_{s\in
 Q'/V}\vr_V((e_0+e_1)^s)=\sum_{V\le Q'}\muqs(V)\sum_{f\in
 F^V}\prod_{s\in
 Q'/V}\vr_V(e_{f(s)}^s)=\sum_{(V,f)\in\fF}\ti\mu(V,f)$$
 where $\ti\mu(V,f)=\muqs(V)\prod_{s\in Q'/V}\vr_V(e_{f(s)}^s)$\,.

 \mn Since $f\in F^V$ implies $fq'\in F^{V^{q'}}\,,\ \fF$ becomes a
 $Q'$-set by $(V,f)^{q'}=(V^{q'},fq')$\,, and we obtain that
 $$\ti\mu((V,f))^{q'}=\ti\mu(((V,f)^{q'})$$
 since, by (2.3), $\ti\mu(V,f)^{q'}=\muqs(V)\prod_{s\in
 Q'/V}\vr_{V^{q'}}(e_{f(s)}^{sq'})=\muqs(V^{q'})\prod_{s_1\in
 Q'/V^{q'}}\vr_{V^{q'}}(e_{(fq')(s_1)}^{s_1})=\ti\mu(V^{q'},fq')\,,$ as $s_1=sq'$ has $(fq')(s_1)=f(s)$.

 \mn We have thus reduced the claimed congruence to
 $$\sum_{(V,f)\in\fF}\ti\mu(V,f)\equiv\ti\mu(V,\textbf{0})+\ti\mu(V,\textbf{1})\mod\tr_{Q'}(\lwa)$$
 where $\fF=\fF(Q')$ and $\textbf{0},\textbf{1}$
 denote the obvious constant functions.

 \mn It now suffices to analyze, for a fixed $f\in F$, the
 $Q'$-orbit sums over $(V,f)\in\fF$. Set $W=\stab_{Q'}(f)$. Note that $W=Q'$ occurs
 only for $f\in F^{Q'}=\{\textbf{0},\textbf{1}\}$, and then in the same way on both sides. Thus we may assume $W<Q'$ from now on.

 \mn Set $e_\ast=\prod_{x\in Q'/W}e_{f(x)}^x\in\lwa$ for a fixed choice of
 coset representatives $x$ of $W$ in $Q'$.
 The $V\le Q'$ for which $(V,f)$ are in $\fF$ are those with
 $f\in F^V$, i.e., $V\le W$. Observe that $(V_1,f)$ and $(V_2,f)$
are in the same $Q'$-orbit if, and only if, $V_1$ and $V_2$ are
conjugate subgroups of $W$. So our sum of $Q'$-orbits involving
$(V,f)$ is \footnote{$\tr_{Q'/V}(e)=\sum_{s\in
Q'/V}e^{s\me}\,,\,V\le Q'\,,\,e\in (\lwa)^V$}
$${\sum_{V\preceq W}}\tr_{Q'/\stab_{Q'}(V,f)}(\ti\mu(V,f))$$
where $\preceq$ means that $V$ runs through a set of
representatives of conjugacy classes of subgroups of $W$. Since
$\stab_{Q'}(V,f)=\stab_{N_{Q'}(V)}(f)=W\cap N_{Q'}(V)=N_W(V)\,,$
it follows that

\mn(2.4)\hsp{15}${{\sum\limits_{V\preceq
W}}}\tr_{Q'/\stab_{Q'}(V,f)}(\ti\mu(V,f))=\tr_{Q'/W}({{\sum\limits_{V\preceq
W}}}\tr_{W/N_W(V)}(\ti\mu(V,f))\,.$

\mn We next analyze $\ti\mu(V,f)$ for each $V\le W$. Decomposing
$W$ as $W=\dot{\bigcup}_yyV$\,, hence
$Q'=\dot{\bigcup}_{x,y}xyV$\,, we obtain
$$\prod_{s\in
Q'/V}e_{f(s)}^s=\prod_{x,y}e_{f(xy)}^{xy}\,\,\dot=\,\,\prod_y(\prod_xe_{f(x)}^x)^y=\prod_ye_\ast^y$$
with $\dot=$ due to $fy\me=f$. Thus
$$\ti\mu(V,f)=\muqs(V)\prod_{s\in
Q'/V}\vr_V(e_{f(s)}^s)=\mu_W(V)\prod_{s\in
W/V}\vr_V(e_\ast^s)\,,$$
hence $$\barr{l}\tr_{W/N_W(V)}\ti\mu(V,f)=\sum_{t\in W/N_W(V)}\ti\mu(V,f)^{t\me}\,\,\dot=\,\,\sum_{t\in W/N_W(V)}\ti\mu(V^{t\me},f)\\
=\sum_{t\in W/N_W(V)}\mu_W(V^{t\me})\prod_{s\in W/V^{t\me}}\vr_{V^{t\me}}(e_\ast^s)\earr$$ with $\dot=$ by $ft\me=f$.

\mn We substitute into our sum (2.4) to get
$$\barr{l}\tr_{Q'/W}(\sum_{V\preceq
W}\tr_{W/N_W(V)}(\ti\mu(V,f))\\
=\tr_{Q'/W}(\sum_{V\le W}\mu_W(V)\prod_{s\in
W/V}\vr_V(e_\ast^s))\in\tr_{Q'/W}(\tr_W(\lwa))=\tr_{Q'}(\lwa)\,,\earr$$
by the induction hypothesis as $W\neq Q'$. Thus, Lemma {\sc b} is
verified.\qed

\bn We now turn to the proof of Proposition 5.

\mn In terms of the map $\vr_V:\lwa\to\lwa$ (where $V=U/A$ and
$A\le U\le G$), Lemma {\sc a} becomes
$$\ver_U^A\res_G^U\ep=\sum_{\pi\in\Si^V}\ti r(\pi)\prod_{s\in Q/V}\vr_V(e_{\pi(s)s\me}^s)$$
with $\ti r(\pi)=\sign(\pi)\prod_{q\in
Q}a_{\pi(q)q\me,q}\in\lwa$\,. Multiplying this by $\muq(V)$ and summing over $V\le Q$ we obtain
\nf{2.5}{\sum_{V\le
Q}\muq(V)\ver_U^A\res_G^U\ep=\sum_{\pi\in\Si}\ti
r(\pi)\sum_{V\le \stab_Q(\pi)}\muq(V)\prod_{s\in Q/V}\vr_V(e_{\pi(s)s\me}^s)}
because  $[\,\pi\in\Si^V\iff V\le\stab_Q(\pi)\,]$\,. We consider the action of conjugation of $Q$ on this sum, starting with the

\mn{\sc Claim\,:}\quad$\ti r(\pi)^q=\ti r(\pi^q)$ \ {\em for all}
$q\in Q$\,. \sn\proof First, $\sign(\pi^q)=\sign(\pi)$ holds since $\{x_i:i\mod b\}$ a cycle of $\pi$ implies $\{x_iq:i\mod b\}$ is a cycle of $\pi^q$. Second,
$$\barr{l}(\prod_{q_1}a_{\pi(q_1)q_1\me,q_1})^q=\prod_{q_1}
a_{\pi(q_1)q_1\me,q_1}^q\sr{1}{=}\prod_{q_1}a_{\pi(q_1)q_1\me,q_1q}
\prod_{q_1}a_{\pi(q_1),q}\me\prod_{q_1}a_{q_1,q}\\
\sr{2}{=}\prod_{q_1}a_{\pi(q_1)q_1\me,q_1q}\sr{3}{=}\prod_{q_1}a_{\pi(q_1q\me)qq_1\me,q_1}=\prod_{q_1}a_{\pi^q(q_1)q_1\me,q_1}\earr$$
with $\sr{1}{=}$ due to the cocycle relation, $\sr{2}{=}$ to $\pi$
permuting the $q_1$ and $\sr{3}{=}$ to the substitution
$q_1\rightsquigarrow q_1q\me$\,.\qed

\mn Continuing with the proof of Proposition 5, the right hand
side of (2.5) is in $\tr_Q(\lwa)$ if
\nf{2.6}{\sum_{V\le\stab_Q(\pi)}\muq(V)\prod_{s\in
Q/V}\vr_V(e_{\pi(s)s\me}^s)\equiv0\mod\tr_{\stab_Q(\pi)}(\lwa)}
holds for all $\pi\in\Si$. Namely, assuming (2.6), its left side
can be written as $\tr_{\stab_Q(\pi)}(\al)$ for some $\al\in\lwa$.
Since $\ti r(\pi)\in (\lwa)^{\stab_Q(\pi)}$\,, it follows that \
$\tr_{Q/\stab_Q(\pi)}(\ti r(\pi)\tr_{\stab_Q(\pi)}(\al))=\tr_Q(\ti
r(\pi)\al)$ \ is the orbit sum of $\pi$ in (2.5), by the above
{\sc Claim}.

\sn We next observe that (2.6) is a consequence of Lemma {\sc b}. To see this, let $\pi\in\Si^V$ be
given, set $Q'=\stab_Q(\pi)$ and $e=\prod_{x\in
Q/Q'}e_{\pi(x)x\me}^x\in \lwa$\,, where $Q=\dot{\bigcup}_xxQ'$\,. Setting
$Q'=\dot{\bigcup}_yyV$\,, then $Q=\dot{\bigcup}_{x,y}xyV$ and the
$V$-term in (2.6) is
$$\muq(V)\prod_{x,y}\vr_V(e_{\pi(xy)(xy)\me}^{xy})\,\,\dot=\,\,\muq(V)\prod_{x,y}\vr_V(e_{\pi(x)x\me}^{xy})=\muqs(V)\prod_{y\in Q'/V}\vr_V(e^y)$$
where $\dot=$ results from $\pi^y=\pi$.

\mn Collecting everything so far, we see that Proposition 5
follows from (2.6), and that this holds because of Lemma {\sc b}.\qed

 \Section{3}{Deriving Proposition 6 from [RW9]}
 For Proposition 6 we
 identify the field $L$ in [RW9] with the fixed field $K^A$ of the given abelian normal open subgroup $A$ in
 $G=G(K/k)$; so $L$ is a Galois extension of $k$ with group $Q$ and the intermediate fields $f$ in $L/k$ correspond to the subgroups $U/A$ of $Q$
 \footnote{our
 $k$ is thus the field $K$ of [loc.cit.]}. Set
 $G_{S'}=G(k_{S'}/k)\,,\,H_{S'}=G(L_{S'}/L)$ \footnote{If Leopoldt's conjecture fails,
 then these groups are not of our type, though abelian.}. Note that the cyclotomic
 $\zl$-extension $L_\infty$ of $L$ is contained in $L_{S'}^+$ and
 abbreviate $G(L_\infty/L)$ by $\Ga_L$.

\def\sp{{_{S'}}}
\def\zlhsp{{\zl[[H\sp]]}}
\def\hpsp{{H_{S'}^+}}
\def\Hsp{{H\sp}}
\mn Let $\la$ be a pseudomeasure on $H\sp$ in the language of
[Se], i.e., an element of the total ring of fractions of the
commutative ring $\zl[[H\sp]]$ so that $(1-h)\la$ is in $\zlhsp$
for all $h\in H\sp$.

\Lemma{c} There is a unique pseudomeasure $\la^A$ on $A$ so that
$(1-a)\la^A$ is the image, under $H\sp\onto A$, of $(1-h)\la$ for
every $h\in H\sp$ with image $a$. Moreover, $\la^A$ is in
$\lba$.\Stop
 We call $\la^A$ the deflation $\defl_{H\sp}^A(\la)$ of the
 pseudomeasure $\la$. Note that if $\la$ is Serre's pseudomeasure
 $\la_L$, then $\la_L^A=\la_{K/L}$.
 \sn\proof Pick an element $h\in H\sp$ with non-trivial image
 under $H\sp\to\Ga_L$. Let $\Ga_h$ be the subgroup of $H\sp$
 topologically generated by $h$ and $M_h\me\zlhsp$ the
 localization of $\zlhsp$ by inverting the multiplicative set
 $M_h=\zl[[\Ga_h]]\setminus l\zl[\Ga_h]]$. Then $\la\in
 M_h\me\zlhsp$ because $(1-h)\la\in\zlhsp$. The natural map
 $\defl_\Hsp^A:\zlhsp\onto\zl[[A]]=\La A$ induces
 $\defl_\Hsp^A:M_h\me\zlhsp\sr{\eta_h}{\to}\lba$. Note that $\lba$ is
 independent, as a subring of $\qq A$, of the choice of $h$, by
 Remark B.

 \sn Then $\eta_h(\la)\in\lba$ is independent of the choice of $h$
 as above, for if $h'$ is another then, with $a$ the image of $h$
 under $H\sp\onto A$,
 $(1-a)\eta_h(\la)=\eta_h((1-h)\la)\,\,\dot=\,\,\eta_{h'}((1-h)\la)=(1-a)\eta_{h'}(\la)\,,$
 with $\dot=$ due to $(1-a)\la\in\zlhsp$. Thus
 $\eta_h(\la)=\eta_{h'}(\la)$ because $1-a\in(\lba)\mal$.
 This common image is $\la^A$.
 \sn Finally, $\la^A$ is a pseudomeaure on $A$ for, if $a\in A$
 then, choosing a preimage $h\in H\sp$ for $a$, we have
 $(1-h)\la\in\zlhsp$ mapping to $(1-a)\la^A$ in $\La A$. Thus, Lemma {\sc c} is shown.\qed

 \mn The Theorem in [RW9] says that there exists a $g\in G\sp$, with
 $g_f$ having image $\neq1$ in $\Ga_f$ for $k\subseteq f\subseteq L$\,, so that
 \nf{3.1}{\sum_{k\subseteq f\subseteq
 L}\muq(G(L/f))\ver_f^L(\ti\la_{g_f})}
 has image in $\tr_Q(\zl[[\hpsp]]$, under $\defl_\Hsp^\hpsp:\zlhsp\to\zl[[\hpsp]]$.

 \mn We first remove the $\ti{\pht{x}}$ on $\la_{g_f}$ in (3.1). Recall that
 $\ti\la_f=2^{-[f:\Q]}\la_f$ implies $\ti\la_{g_f}=2^{-[f:\Q]}\la_{g_f}$. So it suffices to show
 $$\mu_Q(G(L/f))2^{-[f:\Q]}\equiv\mu_Q(G(L/f))2^{-[L:\Q]}\mod|Q|\,,$$
 and this congruence in turn is a consequence of $|Q|$ being
 a power of $l\neq2$ and
 $$\barr{l}\mu_Q(G(L/f))(2^{-[f:\Q]}-2^{-[L:\Q]})=2^{-[L:\Q]}\mu_Q(G(L/f))\Big((2^{[L:f]-1})^{[f:\Q]}-1\Big)\,,\\
 2^{[L:f]-1}\equiv1 \mod l \quad \mr{by} \ 2^l\equiv2\mod l\,,\
 \mr{so} \quad (2^{[L:f]-1})^{[f:\Q]}\equiv 1\mod l[f:\Q]\,,\earr$$
 so it now suffices to show that $\muq(G(L/f))l[f:\Q]$ is divisible by $|Q|=[L:k]$. Due to
 [HIÖ, Corollary 3.5], $\muq(G(L/f))l$ is divisible by $[L:f]$, whence $\muq(G(L/f))l[f:\Q]$ by $[L:f][f:\Q]=[L:\Q]$.

 \mn Next, setting $h=\ver_k^Lg=\ver_f^L\ver_k^fg=\ver_f^Lg_f$\,, we have $(1-g_f)\la_f=\la_{g_f}\in\zl[[G(f\sp/f)]]$ implying that
 $$\la_f\in M_{g_f}\me\zl[[G(f\sp/f)]]$$
 (in the notation of the proof of Lemma {\sc c}) hence
 $$\ver_f^L(\la_f)\in M_h\me\zlhsp$$ because $\ver_f^L=\ver_{G(L\sp/f)^\ab}^\Hsp=\ver_{G(f\sp/f)}^\Hsp$\,.
 Thus $$\sum_{k\subseteq f\subseteq
 L}\muq(G(L/f))\ver_f^L((1-g_f)\la_f)=(1-h)\sum_{k\subseteq f\subseteq
 L}\muq(G(L/f))\ver_f^L(\la_f)$$ implies that every term of the second sum is in $M_h\me\zlhsp$\,. Applying $\defl_\Hsp^A=\defl_\hpsp^A\defl_\Hsp^\hpsp$ to (3.1) then yields that
 \nf{3.2}{\sum_{k\subseteq f\subseteq
 L}\muq(G(L/f))\defl_\Hsp^A\ver_f^L(\la_f)}
 is in $\defl_\hpsp^A(\tr_Q(M_{h^+}\me\zlhsp))=\tr_Q(\lba)$\,, where $h^+=\defl_\Hsp^\hpsp h$\,.

 \mn In order to derive Proposition 6 from this we set $U=G(K/f)$
 and show that
 \nf{3.3}{\defl_\Hsp^A\ver_f^L=\defl_\Hsp^A\ver_{G(f\sp/f)}^\Hsp=\ver_\uab^A\defl_{G(f\sp/f)}^\uab\,.}
 Here, the second equality follows from the commutative diagram
 $$\barr{ccccc}G(L_{S'}/K)&\into&G(L_{S'}/f)&\onto&G(K/f)\\ \pl&&\cup&&\cup\\ G(L_{S'}/K)&\into&G(L_{S'}/L)&\onto&G(K/L)\earr$$
 allowing transfers to be computed by using corresponding coset representatives, giving
 $$\defl_{G(L_{S'}/L)}^{G(K/L)}\ver_{G(L_{S'}/f)}^{G(L_{S'}/L)}=\ver_{G(K/f)}^{G(K/L)}\defl_{G(L_{S'}/f)}^{G(K/f)}\,.$$
 Factoring the transfers through the abelianisations
 $$\barr{l}\defl_{G(L_{S'}/L)}^{G(K/L)}\ver_{G(L_{S'}/f)^\ab}^{G(L_{S'}/L)}\defl_{G(L_{S'}/f)}^{G(L_{S'}/f)^\ab}=\ver_{G(K/f)^\ab}^{G(K/L)}
 \defl_{G(K/f)}^{G(K/f)^\ab}\defl_{G(L_{S'}/f)}^{G(K/f)}\\ =\ver_{G(K/f)^\ab}^{G(K/L)}
 \defl_{G(f_{S'}/f)}^{G(K/f)^\ab}\defl_{G(L_{S'}/f)}^{G(L_{S'}/f)^\ab}\earr$$
 and cancelling the (surjective) maps on the right gives (3.3).

 \sn Combining (3.2) with the commutative square
 $$\barr{cccc}M_{g_f}\me\zl[[G(f\sp/f)]]&\sr{\ver_{G(f\sp/f)}^\Hsp}{\to}&M_{h}\me\zlhsp&\\\daz{\defl_{G(f\sp/f)}^\uab}&&\daz{\defl_\Hsp^A}&\\
 \La_\bullet \uab&\sr{\ver_\uab^A}{\to}&\lba&,\earr$$
 which follows from (3.3), we get the statement of the proposition, by applying Lemma {\sc c} and $\uab=G(K^{[U,U]}/K^U)$. \qed


 \Section{4}{Proof of Theorem 7}
 We fix some notation.
 $A$ is an abelian normal open subgroup of $G=G(K/k)=\group{z}\times G[l]$ containing $\group{z}$ with factor group $Q=G/A$,
 and $\Ga\,,\,C$ are central subgroups of $G$ contained in $A$,
 with $\Ga\simeq\zl$ open and $C$ of exponent $l$.  Further, $U\le G$ is open and contains $C$ and $\group{z}$.

 \mn  Recall that $\Res_G^U$ satisfies diagram
 (D2) and is discussed in \S6 and [RW8, Proposition
 A]. In particular, we know that `Res' is additive, transitive and
 that it
 preserves integrality.
 Also, if $g\in G$, then conjugation by $g$ canonically induces maps
 $\lw U\to\lw U^g$ and $T(\lw U)\to T(\lw U^g)$\,, the latter by $\tau_U(u)^g=\tau_{U^g}(u^g)$ for $u\in U$ \footnote{$\tau_U:\lw U\to T(\lw U)$ is the natural map.}.

 \Lemma{d}\ben\item If $[G:U]=l$ and $\al\in\lw\Ga\,,\,g\in G$\, we have
 $$\Res_{G}^{U}(\al\tau_{G}(g))=\left\{\barr{ll}\al\sum_{x\in G/U}\tau_{U}(g^{x}) &\mr{if}\
g\in U\\ \Psi(\al)\tau_{U}(g^l)&\mr{if}\ g\notin
U\,,\earr\right.$$ where $\Psi:\lw\Ga\to\lw\Ga$ is the continuous
$\zl$-linear ring homomorphism induced by $\ga\mapsto\ga^l$ for
$\ga\in\Ga$. \item For $g\in G\,,\,c\in C$ and $\al\in\lw\Ga$\,,
$$\Res_G^U(\al\tau_G(g(c-1)))=\al\sum_{x\in G/U}\dot\tau_U(g^x(c-1))$$
with $x\in G/U$ meaning $G=\dot{\bigcup}_xxU$ and
$\dot\tau_U:\lwg\to T(\lw U)$ defined by extending $\tau_U$ to
take the value 0 outside of $\lw U$\,.
 \item Let $V\le U$ be
open subgroups of $G$ containing $\group{z}$ and let $g\in
G\,,\,w\in T(\lwg)$\,. Then $(\Res_U^V
w)^g=\Res_{U^g}^{V^g}(w^g)$\,. \een\Stop \proof For $\al\in\Ga$
the claimed formula 1.~is already given at the bottom of
[loc.cit.,p.127]; additivity and continuity then implies it in
general.

\mn For 2., we first remark that $\sum_{x\in
G/U}\dot\tau(g^x(c-1))$ is independent of a special choice of
coset representatives $x$ of $U$ in $G$ and on replacing $g$ by
$g^s$ for $s\in G$. We proceed by induction on $[G:U]$. If $U<G$,
choose $U\le G'< G\,,\,[G:G']=l$, so $G'\lhd G$. Then
$$\barr{l}\Res_G^U(\al\tau_G(g(c-1)))=\Res_{G'}^U\Res_G^{G'}(\al\tau_G(g(c-1)))\\ =\Res_{G'}^U\Big(\left\{\barr{ll}\al\sum_{x''\in G/G'}\tau_{G'}(g^{x''}(c-1)) & \mr{if}\ g\in G'\\ \Psi(\al)\tau_{G'}((g^{x''}c)^l-(g^{x''})^l)&\mr{if} g\notin G'\earr\right.\Big)\\
=\left\{\barr{ll}\sum_{x''\in
G/G'}\Res_{G'}^U\al\tau_{G'}(g^{x''}(c-1))&\mr{if}\ g\in G'\\
0&\mr{if}\ g\notin G'\ .\earr\right.\earr $$ In the first case,
the induction hypothesis turns this into $$\al\sum_{x''\in
G/G'}\sum_{x'\in G'/U}\dot\tau_U(g^{x''x'}(c-1))=\al\sum_{x\in
G/U}\dot\tau_U(g^x(c-1))\,,$$ as required. In the second case,
when $g\notin G'$, we must show that the right hand side of the
assertion is zero, which however is clear since $g^x\in U$ implies
$g^x\in G'$, hence $g\in G'$. This finishes the proof of 2.

 \mn For 3., choose a sequence of
 subgroups $V=V_0<V_1<\cdots<V_n=U$ with $[V_{i+1}:V_i]=l\,,\,0\le i\le n-1$, and then combine transitivity of `Res' with
 induction on $n$ to arrive at
 $$(\Res_U^Vw)^g=(\Res_{V_1}^V(\Res_U^{V_1}w))^g\sr{!}{=}\Res_{V_1^g}^{V^g }((\Res_U^{V_1}w)^g)=\Res_{V_1^g}^{V^g}(\Res_{U^g}^{V_1^g}w^g)=\Res_{U^g}^{V^g}(w^g)\,,$$
 in which the equality $\sr{!}{=}$ still needs to be verified.
 However, here we are in the index $l$ case, so $V^g\lhd V_1^g$, and we can apply 1.

\mn Lemma {\sc d} is established.\qed

 \Lemma{e}  Denote the map $G\to G/C$ by $\ol{\pht{x}}$ and define $\fa\,,\,\fb$ by
 $$\fa\into\lwg\onto\lw\gq\quad\mr{respectively}\quad\fb\into\lw A\onto\lw\aq\ .$$ Then \ben\item[i.]$\fb^Q/\tr_Q\fb\to(\lwa)^Q/\tr_Q(\lwa)$ is injective,
 \item[ii.] $\tau_G(\fa)\into T(\lwg)\onto T(\lw\gq)$ is exact and $Res_G^A\tau_G(\fa)=\tr_Q\fb$\,.\een\Stop
 \proof Applying Tate-cohomology to the sequence defining $\fb$, we see that
 the claimed injectivity i.~is a
 consequence of $H\me(Q,\lw\aq)=0$. In order
 to see this vanishing, choose a central open $\Ga\simeq\zl$ of $\gq$ contained in $\aq$ and pick $Q$-orbit
 representatives $a_i\Ga$ of $\aq/\Ga$. Set $Q_i=\st_Q(a_i\Ga)$.
 For $q_i\in Q_i$, $a_i^{q_i-1}\in\Ga$ has finite order, hence $a_i^{q_i}=a_i$, i.e., $Q_i=\st_Q(a_i)$.
 Hence we have a set of representatives of $\Ga$ in $\aq$
 consisting of $Q$-orbits $Q_i\backslash Q$ for some $Q_i\le Q$ and
 consequently, [RW3, Lemma 5],
 $$H\me(Q,\lw\aq)=\bigoplus_iH\me(Q,\ind_{Q_i}^Q(\lw\Ga))=\bigoplus_iH\me(Q_i,\lw\Ga)=0\,,$$
 as $\lw\Ga$ has $\zl$-torsion $=0$.

 \mn For the first claim in ii.~we only need to check exactness at the middle, or, more precisely, that $\defl_G^\gq \tau_G(v)=0$ implies $\tau_G(v)\in\tau_G(\fa)$. Now, $\defl_G^\gq(v)=\sum_i[\ol w_i,\ol x_i]$ with $w_i, x_i\in\lwg$ implies that $v-\sum_i[w_i,x_i]\in\ker\defl_G^\gq=\fa$\,, so apply $\tau_G$ and arrive at $\tau_G(v)\in\tau_G(\fa)$.
  \mn Regarding the second claim of ii., the elements of $\fa$ are $\lw \Ga$-linear combinations of $\tau_G(g(c-1))$. By 2.~of Lemma {\sc d} $\Res_G^A$ of this equals $\sum_{x\in G/A}\dot\tau_A(g^x(c-1))=\tr_Q(g(c-1))$ if $g\in A$ and 0 if $g\notin A$; note that the $g(c-1)\,,\,g\in A\,,\,c\in C$ generate $\fb$ as $\lw\Ga$-module.

 \sn This finishes the proof of the lemma.\qed
 \Lemma{f} Notation as in Lemma {\sc e}, we have $\Ll(1+\tr_Q\fb)\subset\tr_Q(\lwa)$\,.\Stop
 \proof Let $\be=\sum_{a\in A,c\in
 C}\be_{a,c}a(c-1)$ be an element in $\fb$, where
 $\be_{a,c}\in\lw\Ga$ for some central open $\Ga\simeq\zl$ contained in $A$. Now
 $$\Ll(1+\tr_Q\be)=\fra1l\log\fra{(1+\tr_Q\be)^l}{\Psi(1+\tr_Q\be)}\,,$$
 by the argument given in [RW3,pp.39/40], which also works in the situation when the unit $u$ to which $\Ll$ is applied (see [loc.cit.,p.39,($\ast$)]) is in $\lw A$ rather than in $\La\Ga$ (the ring $\fO$ there is $\zl$ here, so the Frobenius automorphism Fr is trivial). The point is that in [loc.cit.,p.40,($\ast\ast$)] we are on the $A$-level and so we still need only consider degree 1 characters.

 \sn So $\Ll(1+\tr_Q\be)\in \tr_Q(\lwa)$ if $\fra{(1+\tr_Q\be)^l}{\Psi(1+\tr_Q\be)}\equiv1\mod l\tr_Q(\lwa)$.
 Since $(1+\tr_Q\be)^l\equiv 1+(\tr_Q\be)^l\mod l\tr_Q(\lw A)$, it suffices to show that $(\tr_Q\be)^l\equiv\Psi(\tr_Q\be)\mod
 l\tr_Q(\lwa)$\,.

 \sn Now, as $\Psi(a)=a^l$ for $a\in A$ [loc.cit.,p.33], $\Psi(\tr_Q\be)=\sum_{a,c}\Psi(\be_{a,c})(\tr_Q((ac)^l)-\tr_Q(a^l))=0$\,,
 since $(ac)^l=a^lc^l=a^l$, and we are left to check that $(\tr_Q\be)^l\equiv0\mod
 l\tr_Q(\lwa)$\,. But $\tr_Q\be=\sum_c\ka_c(c-1)$ with
 $\ka_c=\tr_Q(\sum_a\be_{a,c}a)\in\tr_Q(\lwa)$\,, so
 $$(\tr_Q\be)^l\equiv\sum_c\ka_c^l(c-1)^l\mod l\tr_Q(\lwa)$$
 as $\tr_Q(\lwa)$ is an ideal in $(\lwa)^Q$ and
 $c\in(\lwa)^Q$\,. Thus, $(c-1)^l\equiv0\mod
 l(\lwa)^Q$ establishes Lemma {\sc f}.\qed

\bn We are now in a position to prove Theorem 7.

\mn If $U$ is an open subgroup of $G$ containing $\group{z}$ and
$N$ a finite normal subgroup of $U$, we write $t_{U/N}$ for
$t_{K^N/K^U}$; similarly, we write $\la_A$ for $\la_{K/K^A}\in(\lw
A)\mal=K_1(\lw A)$. We recall from [RW8, Lemma 2] and [RW5, Lemma
2.1] that $\Res_G^U t_G=t_U\,,\,\defl_U^{U/N}t_U=t_{U/N}$\,, in
this notation.

 \sn As above, denote going modulo
 $C$ by $\ol{\pht{x}}$\,, so $\ol G=G/C$. Also recall the short exact sequences
 $$\fa\into\lwg\onto\lw\gq\quad,\quad\fb\into\lw A\onto\lw\aq\,.$$

 \mn For the definitions to follow we use the commutative square

 \mn\bmp{6cm}$\barr{ccc}(\lwg)\mal&\onto&(\lw\gq)\mal\\ \sda&&\sda\\ K_1(\lwg)&\to&K_1(\lw\gq)\earr$\emp
 \hspace{10mm}\bmp{7cm} in which the (natural) vertical maps are surjective,
 since $\lwg\,,\,\lw\gq$ are semi-local rings. Moreover, the top horizontal map is
 surjective as well, because $\ker(\lwg\to\lw\gq)\subset\mr{rad}(\lwg)$\,.\emp

 \bn By Proposition 3, $t_\gq\in T(\lw\gq)$ implies $L_{K^C/k}=\Det\,\theta$
 with $\theta\in K_1(\lw\gq)$. Observe that $\res_\gq^\aq\theta=\la_\aq$\,, because
 $$\Det(\res_\gq^\aq\theta)=\res_\gq^\aq(\Det\theta)=\res_\gq^\aq(L_{K^C/k})=L_{K^C/K^A}=\Det\la_\aq$$
 and $SK_1(\lw\aq)=1$, as $\lw\aq$ is commutative semilocal [CR, (45.12)].

 \mn The above square gives a $\vt\in K_1(\lwg)$ with $\defl_G^\gq\vt=\theta$. Define
 $\vt'\in K_1(\lw A)$ by $\res_G^A\vt=\la_A\vt'$\,, hence $\defl_A^\aq\vt'=1$.

 \sn Further define\quad $\xi=\Ll(\vt)\,,\,\xi'=\Ll(\vt')$\,. Then, using diagrams (D1) and (D2),
 $\xi\in T(\lwg)$ has \ $\defl_G^\gq\xi=t_\gq\,,\,\Res_G^A\xi=t_A+\xi'$ and $\defl_A^\aq\xi'=0$.

 \mn The exact sequences displayed above give rise to the commutative diagram
 $$\barr{ccccc}\tau_G(\fa)&\into&T(\lwg)&\onto&T(\lw\gq)\\
 \da&&\daz{\Res_G^A}&&\daz{\Res_\gq^\aq}\\
 \fb&\into&\lw A&\onto&\lw\aq\earr$$ with top sequence exact by Lemma {\sc e}.
 We need to modify $\xi$ by adding an element $\al\in\tau_G(\fa)$ (so without changing $\defl_G^\gq\xi$) to arrange that $\Res_G^A(\xi+\al)=t_A$\,, i.e., we need to prove that $\xi'\in \Res_G^A(\tau_G(\fa))$.

 \mn Now, $\la_{A}$ is $Q$-invariant, by the proof of [RW5, Lemma
 3.1], and $\res_G^A$ takes $Q$-invariant values, whence $\vt'\in 1+\fb^Q.$
 We claim that
 \nf{4.1}{\vt'\in 1+\tr_Q\fb\,.}
 Write $\defl_U^\uab\res_G^U\vt=\la_\uab \vt_U'$ for $A\le U \le G$.
 Then $\vt_U'\in K_1(\lw\uab)=(\lw\uab)\mal$ has
 $\defl_\uab^{U/\uuc\cdot C}\vt_U'=1$\, because
 $$\barr{l}\la_{U/[U,U]C}\defl_\uab^{U/[U,U]C}(\vt_U')=\defl_\uab^{U/[U,U]C}\defl_U^\uab\res_G^U\vt
 =\defl_\uq^{U/[U,U]C}\defl_U^\uq\res_G^U\vt\\ =\defl_\uq^{U/[U,U]C}\res_\gq^\uq\defl_G^\gq\vt
 =\defl_\uq^{U/[U,U]C}\res_\gq^\uq\theta=\la_{U/[U,U]C}\,,
 \earr$$ with the last equality by $SK_1(\lw U/[U,U]C)=1$, as before.

 \sn Summing up,
 $\vt_U'\in\ker(\,(\lw\uab)\mal\to(\lw U^{U/\uuc\cdot C})\mal\,)$
 so $\vt_U'=1+\al$ where $\al$ is a $\lw\Ga$-linear combination of elements
 $\ti u(\ti c-1)$ with $u\in U\,,\,c\in C$ and
 $\ti{\pht{x}}:U\to\uab$ the canonical map (with $\Ga\simeq\zl$ some central open subgroup of $G$ contained in $A$). Then, if $A<U$,
  $\ver_\uab^A$ takes $c$ to $c^{[U:A]}=1$ and thus $\ti
 u(\ti c-1)$ to zero. As a first result we
 therefore have $$U\neq A\implies\ver_U^A\res_G^U \vt=\la_\uab\ .$$
 Now insert $\vt$ into the ``Möbius-Wall'' congruence of Proposition 5 and obtain
 $$\la_A\vt'+\sum_{A<U\le
 G}\mu_Q(U/A)\ver_\uab^A\la_\uab\equiv0\mod\tr_Q(\lwa)\ .$$
 Comparing this \footnote{It is only here that Propositions 5 and 6 make their appearance.} with the abelian pseudomeasure congruence
 $$\sum_{A\le U\le G}\muq(U/A)\ver_\uab^A(\la_\uab)\equiv0\mod \tr_Q(\lw A) $$
 of Proposition 6 gives \
 $\vt'\equiv1\mod\tr_Q(\lwa)$\,, as $\la_A$ is a unit in $(\lwa)^Q$ and $\tr_Q(\lwa)$ an ideal. Therefore
 $$\vt'\in(1+\fb^Q)\cap(1+\tr_Q(\lwa))=1+(\fb^Q\cap\tr_Q(\lwa))=1+\tr_Q\fb\,,$$
 with the last equality due to Lemma {\sc e}. This proves our claim (4.1).

 \mn  Turning back to the proof of the theorem, we know that $\xi'=\Ll(\vt')$ is in $\fb^Q$. By the claim and Lemma {\sc f}, we also have
 $\xi'\in \tr_Q(\lw A)$ hence $\xi'\in\fb^Q\cap\tr_Q(\lw A)=\tr_Q(\fb)=\Res_G^A(\tau_G(\fa))$\,, by Lemma {\sc e}.

 \sn Thus the proof of Theorem 7 is complete.\qed

 \Section{5}{Proof of the Theorem}
 Recall the notation of the beginning of \S4, so $G=G(K/k)=\group{z}\times G[l]$ is
 $l$-elementary, with $\group{z}$ a finite cyclic group of order prime to $l$ and $G[l]$ a pro-$l$ group,
 and $A$ is an abelian normal open subgroup of $G$ containing $\group{z}$.

 \sn We define $\fc_U^\ab$ by the exact
 sequence $0\to\fc_U^\ab\to\lw(\uab)\to\lw(U/C[U,U])\to0$\,.

 \Lemma{g} Let $C$ have order $l$ and $U\ge A$ satisfy $C\cap[U,U]=1$\,.
 Denote the normalizer of $U$ in $G$ by $N=N_G(U)$ and let $Y$ be a set of representatives of $N/U$-orbits in
 $U/\Ga C[U,U]$\,. Then $\tr_{N/U}(\fc_U^\ab)$ has $\lw\Ga$-basis
 $$\{\tr_{N/U}(\ti y(\ti c-1)):y\in Y_1,\ 1\neq c\in C\}$$
  where $Y_1$ is the set of $y\in Y$ that have preimage $\ti y$ in
  $U/[U,U]$ which is fixed by $\st_{N/U}(y)$\,, and $\ti c$ is
  the image of $c\in C$ in $U/[U,U]$.\Stop
  For the proof, we use $C\cap[U,U]=1=\Ga\cap[U,U]$ to
  identify $C,\Ga$ with their images in $\uab$ (hence $c$ with
  $\ti c$). We investigate the $N/U$-structure of
  $0\to\fc_U^\ab\to\lw\uab\to\lw(\uab/C)\to0$ via the
  $\lw\Ga$-bases coming from the $N/U$-action on
  $C\into\uab/\Ga\onto\uab/\Ga C$ by [RW3, Lemma 5].

  \sn Now $Y$ is a set of representatives of $N/U$-orbits on
  $\uab/\Ga C$. If $\hat y$ is a preimage of $y\in Y$ under
  $\uab/\Ga\to\uab/\Ga C$ then $\st_{N/U}(y)$ either fixes $\hat
  y$ (in case 1) or moves $\hat y$ (in case 2); moreover this case
  distinction is independent of the choice of $\hat y$. This
  permits us to analyze the map $\lw\uab\to\lw(\uab/C)$ one $y\in
  Y$ at a time in terms of the map of $N/U$-sets from the
  preimage of the $N/U$-orbit of $y$ to the $N/U$-orbit of $y$ it
  induces. This is because of the permutation $\lw\Ga$-basis given
  by choosing preimages $\ti y$ of $\hat y$ under
  $\uab\to\uab/\Ga$ with $\st_{N/U}(\ti y)=\st_{N/U}(\hat y)$, as
  in the proof of Lemma {\sc e}.\sn  Thus, in case 1, the preimage of the
  $N/U$-orbit of $y$ is $\dot{\bigcup}_{c\in C}(N/U$-orbit of $\hat
  yc)$\,, so $l$ copies of $\fra{N/U}{\st_{N/U}(y)}$ as $N/U$-sets, and the map is $\hat
  y^nc\mapsto y^n$ for $n\in N,c\in C$. So the kernel on
  $\lw\Ga$-permutation modules has $\lw\Ga$-basis $\{\ti
  y^n(c-1):n\in\fra{N/U}{\st_{N/U}(y)}\,,\,1\neq c\in C\}$\,.

  \sn Similarly, in case 2, the preimage of the $N/U$-orbit of $y$
  is the $N/U$-orbit of $\hat y$\,: here $\hat y^z=\hat
  y\ga_y(z)$,
  with $\ga_y$ a homomorphism $\st_{N/U}(y)\onto C$, has
  $\st_{N/U}(\hat y)$ as its kernel. Now the kernel on
  $\lw\Ga$-permutation modules has $\lw\Ga$-basis $\{\ti y^n-\ti
  y:n\in\fra{N/U}{\st_{N/U}(\hat y)}\}$\,.

  \sn Hence $\fc_U^\ab$ has $\lw\Ga$-basis the union of these over
  $y\in Y$, and $\tr_{N/U}(\fc_U^\ab)$ has the claimed
  $\lw\Ga$-basis since $Y_1$ consists of the $y\in Y$ in case 1 and
  $\tr_{N/U}(\ti y^n-\ti y)=0$ for all $y\in Y$ in case 2. This
  proves the lemma.\qed

\def\vab{{V^\ab}}
 \Lemma{h} If $v\in T(\lwg)$ has $\defl_U^\uab\Res_G^Uv=0$
 for all subgroups $U$ of $G$ containing $A$, then $v=0$.\Stop
 The proof is by induction on $[G:A]$. Fix a central open $\Ga\simeq\zl$ in $A$ and an $n$ so that $l^n\equiv1\mod|z|$ and $l^n$ is an
 exponent of $G[l]/\Ga$. By the diagram of ii) of
 Lemma {\sc i} and $\defl_G^\gab v=0$, there is an $\om\in
 K_1(\lwg)$ so that $\Ll(\om)=v$ and  $\defl_G^\gab\om$ is a torsion element (e.g. 1) of $K_1(\lw\gab)$.

 \sn Consider $\Res_G^Uv$ with $U\ge A$ having index $l$ in $G$. Then $\defl_V^\vab\Res_U^V(\Res_G^Uv)=\defl_V^\vab\Res_G^Vv=0$ so the induction hypothesis yields $\Res_G^Uv=0$. Thus $\Ll(\res_G^U\om)=\Res_G^U(\Ll\om)=0$ implies $\LL(\Det(\res_G^U\om))=\Tr(\Ll(\res_G^U\om))=0$\,, hence $\Det(\res_G^U\om)(\chi_1)^l=\Psi((\Det(\res_G^U\om))(\psi_l\chi_1)$ for all characters $\chi_1$ of $U$, by the definition of $\LL$.

 \sn If $\chi_1$ is an irreducible character of $U$ with kernel containing $\Ga$ then $\chi_1=\be\ot\varpi$ with $\be,\varpi$ irreducible characters of $U$ with kernels containing $U[l],\group{z}\Ga$, respectively. Note that $\psi_l^n\chi_1=\psi_l^n\be\ot\psi_l^n\varpi=\be\ot1=\be$\,, hence $$(\Det(\res_G^U\om)(\chi_1)^{l^n}=\Psi^n(\Det(\res_G^U\om)(\psi_l^n\chi_1))=\Psi^n((\res_G^U(\Det\om))(\be))=\Psi^n((\Det\om)(\ind_U^G\be))$$ is torsion. This holds because $U\lhd G$  $l$-elementary implies that $\be=\res_G^U\be'$ where $\be'=\infl_\gab^G\be''$ hence $\ind_U^G\be=\be'\ind_U^G1$ and $\ind_U^G1=\sum_{i=1}^l\al_i$ with $\al_i$ the irreducible characters of $G$ having $\res_G^U\al_i=1$\,; now $\al_i=\infl_\gab^G\al_i'$ so $$(\Det\om)(\ind_U^G\be)=(\Det\om)(\sum_{i=1}^l\be'\al_i)=
 \prod_{i=1}^l(\Det\om)(\infl_\gab^G(\be''\al_i'))=\prod_{i=1}^l\Det(\defl_G^\gab\om)(\be''\al_i')$$
 is torsion. Thus $(\Det(\res_G^U\om))(\chi_1)$ is torsion for all such $\chi_1$ and $U\ge A$ of index $l$ in $G$.

 \sn Now if $\chi$ is one of the finitely many irreducible characters of $G$ with $\Ga\subseteq\ker(\chi)$ then (by Clifford theory), either $\chi=\ind_U^G\chi_1$ with
 $U\ge A$ of index $l$ when $(\Det\om)(\chi)=(\res_G^U(\Det\om))(\chi_1)$ is torsion, or $\chi=\infl_\gab^G\al$ when $(\Det\om)(\chi)=\Det(\defl_G^\gab\om)(\al)$ is again torsion. Every irreducible character of $G$ has the form $\chi\ot\rho$ with such a $\chi$ and $\rho$ of type W hence $(\Det\om)(\chi\ot\rho)=\rho^\sharp((\Det\om)(\chi))$ torsion of order at most that of $(\Det\om)(\chi)$. Thus $\Det\om$ is a torsion element in $\HOM(R_lG,(\lw\Ga_k)\mal)$ and so $\Tr(v)=\Tr(\Ll \om)=\LL(\Det\om)=0$ implies $v=0$. \qed

 \Remark{C} When $G$ is abelian pro-$l$, Proposition 5.1 of [RW6] gives a description of the kernel of $\Ll$ on $(\lwg)\mal$. This can be (and originally was) used to prove the lemma for pro-$l$ groups. The present proof is shorter for $l$-elementary groups than extending that proposition from $\La_\wedge$ to $\lwb$ (for the notation compare the proof of Lemma {\sc i}).

 \bn We next state the

 \mn{\sc Uniqueness Principle.}\quad{\em If $\xi\in T(\lwg)$ has $\defl_U^\uab\Res_G^U\xi=t_\uab$
 for all subgroups $U$ of $G$ containing a fixed abelian
 normal open subgroup $A$ of $G$, then $\xi=t_G$. In particular, $t_G\in T(\lwg)$.}

 \mn This is because $l^nt_G$ is integral for large enough natural $n$
 (see [RW4, Proposition 9]). Setting $v=l^nt_G-l^n\xi$ we see that
 $\defl_U^\uab\Res_G^Uv=l^nt_\uab-l^nt_\uab=0$\,, so
 $v=l^n(t_G-\xi)=0$ by Lemma {\sc h}. However, $T(\lwg)$ is torsionfree as `Tr' is injective.\qed

 \bn Now, were the {\sc Theorem} false, there would
 exist an extension $K/k$ for which the Galois group $G$ would
 have commutator subgroup $[G,G]$ of {\em minimal} order; among
 these groups $G$ there would be one with centre $Z(G)$ of {\em
 minimal} index $[G:Z(G)]$.

 \sn Since $[G,G]\neq1$, by [RW3, Theorem 9], and $[G,G]$ is an $l$-group as $G$ is $l$-elementary we may choose a central subgroup
 $C\le[G,G]$ of order $l$ in $G$, and then a {\em maximal} abelian
 normal subgroup $A$ of $G$, necessarily containing $C$ and $\group{z}$. We also fix a central open
 $\Ga\simeq\zl$ inside $A$.  Now
 Theorem 7 guarantees the existence of
 \nf{5.1}{\xi\in T(\lwg)\ \mr{with}\ \defl_G^\gq \xi=t_\gq\ \mr{and}\
 \Res_G^A\xi=t_A\,,}
 where, as before, $\ol{\pht{x}}$ denotes going modulo $C$. To defeat the
 counterexample $G$ it suffices, by the {\sc Uniqueness Principle}, to
 find such a $\xi$ so that $\defl_U^\uab\Res_G^U\xi=t_\uab$ for all
 subgroups $U\ge A$ of $G$. Observe that this already holds for
 $U$ with $[U,U]\ge C$\,: for then
 $\defl_U^\uab\Res_G^U\xi=\defl_{U/C}^\uab\defl_U^{U/C}\Res_G^U\xi=\defl_\uq^\uab\Res_\gq^\uq\defl_G^\gq \xi=\defl_\uq^\uab\Res_\gq^\uq
 t_\gq=t_\uab\,.$

 \sn On the other hand, for $U\ge A$ with $[U,U]\not\ge C$, hence
 $C\cap[U,U]=1$, then $|[U,U]|<|[G,G]|$ implies \footnote{recall
 that $t_U=t_{K/K^U}$ by the Galois correspondence}
 $t_U\in T(\lw U)$,
 by our hypothesis on $G$, permitting us to define
 $$\xi_U=\Res_G^U\xi-t_U\in T(\lw U)$$
 and to define the {\em support} of $\xi$ by
 $$\mr{supp}(\xi)=\{U\ge A:C\cap[U,U]=1 \ \mr{and} \ \xi_U\neq0\}\,.$$   To investigate $U\in\supp(\xi)$, we state

 \mn{\sc Claim 1\,:}\quad \ben\item[a)]{\em If $A\le V\le U$ and $C\cap[U,U]=1$\,, then $\Res_U^V\xi_U=\xi_V$ and $\xi_V^g=\xi_{V^g}$ for $g\in G$.}
 \item[b)] {\em $G$ acts on $\supp(\xi)$ by conjugation.}
 \item[c)] $A\notin\supp(\xi)$\,.
 \een

 \mn \proof Recall that $t_U=\Res_G^Ut_G$\,. Now, a) results from
 $$\barr{l} \Res_U^V\xi_U=\Res_U^V(\Res_G^U\xi-\Res_G^Ut_G)=\Res_G^V\xi-\Res_G^Vt_G=\xi_V\quad\mr{and}\\
 \xi_V^g=\Res_G^V(\xi-t_G)^g=\Res_{G^g}^{V^g}(\xi^g-t_G^g)=\Res_G^{V^g}(\xi-t_G)=\xi_{V^g}\ ,\earr$$
 by 3.~of Lemma {\sc d}, which at the same time implies b); c) follows from (5.1) and the definition of `supp'.\qed

 \mn Moreover, we let, as in Lemma {\sc g}, $N=N_G(U)$ be the normalizer
 of $U$ in $G$. Define $\fc_U$ and $\fc_U^\ab$ by the exact
 sequences $$\fc_U\into\lw(N/[U,U])\onto\lw(U/C[U,U])\quad\mr{and}\quad
 \fc_U^\ab\into\lw(\uab)\onto\lw(U/C[U,U])\,.$$

 \sn {\sc Claim 2\,:}\quad {\em If $U\in\supp(\xi)$\,,
 then $\defl_U^\uab \xi_U\in\tr_{N/U}(\fc_U^\ab)$\,.}

 \mn\proof We first note that $t_{N/[U,U]}$ is in
  $T(\lw(N/[U,U]))$\,: for the commutator subgroup $[N,N]/[U,U]$
  of $N/[U,U]$ has smaller order than $[G,G]$ unless $[U,U]=1$ and
  $[N,N]=[G,G]$, in which case $A\notin\supp(\xi)$ implies
  $N<G$, because $A$ is {\em maximal} abelian normal in $G$; but
  then $Z(G)\le N$ implies $[N:Z(N)]\le[N:Z(G)]<[G:Z(G)]$\,,
  contrary to the minimality hypothesis on $G$.

  \sn Writing $\defl_N^{N/[U,U]}\Res_G^N\xi=t_{N/[U,U]}+z_U$\,, with
  $z_U\in T(\lw(N/[U,U]))$\,, and $\tau=\tau_{N/[U,U]}$, we consider the commutative diagram
  $$\barr{cccccc}&&T(\lw N)&\onto&T(\lw\nq)&\\ &&\sda&&\sda&\\
  \tau(\fc_U)&\into&T(\lw(N/[U,U]))&\onto&T(\lw(N/C[U,U]))&,\earr$$
  with all surjective maps deflations and exact bottom row by ii.~of Lemma {\sc e} applied to $N/[U,U]\ge U/[U,U]\ge C[U,U]/[U,U]$
  in place of $G\ge A\ge C$; moreover, we also obtain $\Res_{N/[U,U]}^\uab\tau(\fc_U)=\tr_{N/U}(\fc_U^\ab)$.

  \mn Since $\Res_G^N\xi\in T(\lw N)$ has
  $\defl_N^\nq\Res_G^N\xi=\Res_\gq^\nq\defl_G^\gq \xi=\Res_\gq^\nq
  t_\gq=t_\nq$\,, the diagram implies $z_U\in\tau(\fc_U)$. Thus
  $$\barr{l}\defl_U^\uab
  \xi_U+t_\uab=\defl_U^\uab(\Res_G^U\xi-t_U)+t_\uab=\defl_U^\uab\Res_N^U(\Res_G^N\xi)=\\
  \Res_{N/[U,U]}^\uab\defl_N^{N/[U,U]}(\Res_G^N\xi)=\Res_{N/[U,U]}^\uab(z_U+t_{N/[U,U]})=
  \Res_{N/[U,U]}^\uab z_U+t_\uab\earr$$ implies $\defl_U^\uab
  \xi_U=\Res_{N/[U,U]}^\uab z_U\in\Res_{N/[U,U]}^\uab\tau(\fc_U)$.
  Combining with the previous paragraph yields the claim.\qed

  \mn Now continuing with the proof, it follows that the {\sc Theorem} holds if we can modify $\xi$, subject to
 (5.1) holding, so that supp$(\xi)$ is empty. Since this is not
 possible for our $G$, by hypothesis, there must exist an $\xi$ for
 which supp$(\xi)$ has {\em minimal cardinality} $\neq0$.

 \mn Since $\supp(\xi)$ is non-empty it contains a $U$ with
  {\em minimal} $[U:A]$. By Lemma {\sc g} we may write
  \nf{5.2}{\defl_U^\uab \xi_U=\sum_{y\in Y_1,\,1\neq c\in
  C}\al(y,c)\tr_{N/U}(\ti y(\ti c-1))} with unique $\al(y,c)$ in
  $\lw\Ga$. As $A\notin\supp(\xi)$ every $V$ with $A\le
  V<U$ has $C\cap[V,V]=1$ and $\xi_V=0$; equivalently, we may
  restrict attention to maximal such $V$, hence may assume that
  $V$ has index $l$ in $U$ (so, in particular, $[U,U]\lhd V$). Now
  $$\Res_\uab^{V/[U,U]}(\defl_U^\uab
  \xi_U)=\defl_V^{V/[U,U]}\Res_U^V\xi_U=\defl_V^{V/[U,U]}\xi_V=0\,,$$
  by a) of {\sc Claim 1}. We
  apply $\Res_\uab^{V/[U,U]}$ to (5.2), observing (by 2.~of Lemma {\sc d}) that, for
  $u\in\uab$, $\Res_\uab^{V/[U,U]}(u(\ti c-1))=lu(\ti c-1)$ if
  $u\in V/[U,U]$, and $=0$, if not, because $\uab$ is abelian. It
  follows that $$0=\sum_{y\in Y_1\cap (V/\Ga
  C[U,U]),\,c\neq1}l\al(y,c)\tr_{N/U}(\ti y(\ti c-1))$$
  whenever $A\le V\le U\,,\,[U:V]=l$\,, and,
  therefore, that $\al(y,c)=0$ unless $y\in Y_1$ is not in $V/\Ga
  C[U,U]$\,.

  \sn In particular, if $U/A$ is non-cyclic, then every element of $U/A$ is
  contained in a maximal $V/A$ for some $V$, hence
  $\defl_U^\uab \xi_U=0$. But then the {\sc Uniqueness Principle}, applied
  to to $A\lhd U$ instead of $A\lhd G$, implies that $\xi_U=0$. Thus
  $U\notin\supp(\xi)$, contrary to assumption.

  \sn It follows that our $U$ with minimal $[U:A]$ in $\supp(\xi)$
  has $U/A$ cyclic. Then $[U,U]\le A$\,, $U/A\simeq\fra{U/\Ga C[U,U]}{A/\Ga C[U,U]}$\,, and
  \nf{5.3}{\defl_U^\uab \xi_U=\sum\limits_{y\in Y_1\,,\group{yA}= U/A\,,\,1\neq
  c\in C}\al(y,c)\tr_{N/U}(\ti y(\ti c-1))\,,}
  because $U/A$ now has a unique maximal subgroup $V/A$ and so $y\in Y_1$ is not in $V/\Ga C[U,U]$ precisely when $yA$ generates $U/A$.

  \mn Now our {\sc Theorem} essentially follows from the next result.

  \mn {\sc Claim 3\,:}\quad {\em Assume that $U/A$ is cyclic. Set, in the notation of (5.3),
  $$\xi''=\sum_{y\in Y_1\,,\,\group{yA}= U/A\,,\,1\neq c\in
  C}\al(y,c)\tau_G(y'(c-1))\quad in \quad T(\lwg)\,,$$ with
  preimages $y'\in U$ of $\ti y$ under $\defl_U^\uab$. Then
  \ben\item[i.] $\defl_U^\uab\Res_G^U\xi''=\defl_U^\uab \xi_U$\,, and
  \item[ii.] if $A\le U_1\le G$ then $\Res_G^{U_1}\xi''\neq0$ implies $\exists\,g\in G:U^g\le\ue$\,.\een}

  \mn\proof Recall that the $y$ in the $\xi''$-sum have $\group{yA}=U/A$.  Applying 2.~of Lemma {\sc d} gives
  $$\Res_G^U\xi''=\sum_{y\in Y_1\,,\,\group{yA}= U/A\,,\,1\neq c}\al(y,c)\sum_{x\in G/U}\dot\tau_U((y')^x(c-1))\,.$$
  Note that $(y')^x\in U$ implies $(yA)^x\in U/A$, hence $(U/A)^x=U/A$, i.e., $x\in N$. Now we have
  $$\Res_G^U\xi''=\sum_{y,c}\al(y,c)\sum_{x\in N/U}\tau_U((y')^x(c-1))\,,$$
  hence applying $\defl_U^\uab$ gives i.

  \mn For ii.~note that $y$ still has $\group{yA}=U/A$, but we now apply 2.~of Lemma {\sc d} with $U$ replaced by $\ue$. Some term $\Res_G^\ue\tau_G(y'(c-1))$ in this sum must be $\neq0$, by hypothesis; but this term is $\sum_{x_1\in G/U}\dot\tau_\ue((y')^{x_1}(c-1))$\,, so we must have a non-zero term here, i.e., $(y')^{x_1}\in U_1$ for some $x_1$. Now $(U/A)^{x_1}\subseteq U_1/A$ implies $U^{x_1}\subseteq U_1$.\qed

 \mn We apply {\sc Claim 3} and set $\xi'\df \xi-\xi''$. Then $\Res_G^A\xi'=t_A$\,, by
 ii.~with $U_1=A$\,; moreover,
  due to the appearance of the elements $1\neq c\in C$ in $\xi''$,
  $\defl_G^{G/C}\xi'=t_{G/C}$\,;
  thus $\xi'$ satisfies (5.1). Further, $\supp(\xi')\subseteq\supp(\xi)$\,: for if $\ue\in\supp(\xi')$ then $C\cap[U_1,U_1]=1$ and $\xi_\ue'\neq0$, hence $\xi_\ue=\xi_\ue'+\Res_G^\ue\xi''$ is nonzero unless $\Res_G^\ue\xi''\neq0$\,; but in that case $U_1\supseteq U^g$ for some $g\in G$ by ii., hence {\sc Claim 1} implies $\xi_\ue\neq0$, as $\Res_\ue^{U^g}\xi_\ue=\xi_{U^g}\neq0$ by $U\in\supp(\xi)$.
  But now i.~and ii.~of {\sc Claim 3} imply $\defl_U^\uab\xi_U'=0$ and $\xi_\ue'=0$ for $A\le U_1< U$, hence $\xi_U'=0$ by the {\sc Uniqueness Principle}\,. Thus  $U\notin\supp(\xi')$, which contradicts
  the minimal cardinality of $\supp(\xi)$ and therefore finishes the proof of the {\sc Theorem}.\qed

 \Section{6}{$l$-elementary groups}
Recall that $G=G(K/k)=\group{z}\times G[l]$ is $l$-elementary.

\Lemma{i} \ben\item[i)] The logarithm $\Ll:K_1(\lwg)\to T(\qwg)$
of diagram (D1) has image in $T(\lwg)$.\item[ii)] Let $\fa$ be the
kernel of $\defl_G^\gab:\lwg\onto\lw\gab$. Then the commutative
diagram
\mn\bmp{7cm}$\barr{ccccc}1+\fa&\into&(\lwg)\mal&\onto&(\lw\gab)\mal\\ \sda&&\daz{\Ll}&&\daz{\Ll^\ab}\\
\tau_G(\fa)&\into&T(\lwg)&\onto&\lw\gab\earr$\emp\hspace{10mm}\bmp{5cm}
has exact rows and surjective left vertical map.\emp\een \Stop
\proof For i), abbreviate $G[l]$ as $U$. Each
$\ql\clo$-irreducible character $\be$ of $\group{z}$ induces a
$\zl$-algebra homomorphism $\zl[\group{z}]\to\zl\clo$ with image
$\zl[\be]$, hence surjective ring homomorphisms
$$\ql[\group{z}]\onto\ql(\be)\,,\,\lwg\onto\zl[\be]\ot_\zl\lw
U\df\lwb U\,,\,\qwg\onto\qwb U\,.$$ Applying the functors $K_1$
and $T$ gives the southeast and southwest arrows of the diagram
{\footnotesize $$\barr{ccccccc}K_1(\lwg)&&&\lto&&&T(\qwg)\\
&\searrow&&&&\swarrow&\\ &&K_1(\lwb U)&\to&T(\qwb U)&&\\ \dal&&\da&&\da&&\dal\\
&&\HOM^{(\be)}(R_lU,(\lwc\Ga_k)\mal)&\to&\Hom^{(\be)}(R_lU,\qwc\Ga_k)&&\\
&\nearrow&&&&\nwarrow&\\
\HOM(R_lG,(\lwc\Ga_k)\mal)&&&\lto&&&\Hom^\ast(R_lG,\qwc\Ga_k)\earr$$}
\hspace{-1mm}with large square from diagram (D1) of \S2, and small
square [RW3, 2.~of Proposition 11] with unramified coefficients
$\zl[\be]$, which are abbreviated by the superscript $\be$. The
northwest and northeast arrows $f\mapsto f^{\be}$ are defined by
$f^{\be}(\varpi)=f(\be\ot\varpi)$\,.


 \sn To see that the left quadrangle commutes \footnote{this can also be obtained from [RW4, Theorem 1]} let $H'=H\cap U$
(recalling that $H$ is the kernel of $G\to\Ga_k$), hence
 $H=\group{z}\times H'$, and let $\be(x)$ denote the image of
 $x\in\lwg$ in $\lwb U$. We must check that $(\Det
 x)^{\be}(\varpi)=(\Det\be(x))(\varpi)$\,, i.e.~\footnote{compare also [RW2, Proposition 6; RW3, Lemma
 2]},
 ${\det}_{\qwc\Ga_k}(x\mid\fV_{\be\ot\varpi})=\det(\be(x)\mid\fV_\varpi)\,.
 $ Here,
 $\fV_\varpi=\Hom_{\ql\clo[H']}(V_\varpi,\ql\clo\ot_{\ql(\be)}\qwb
 U)=\Hom_{\ql\clo[H']}(V_\varpi,\qwc U)$ and
 $$\barr{l}\fV_{\be\ot\varpi}=\Hom_{\ql\clo[H]}(V_{\be\ot\varpi},\qwc G)=\Hom_{\ql\clo[\group{z}]\ot_{\ql\clo}\ql\clo[H']}(V_\be\ot_{\ql\clo}V_\varpi,
 \ql\clo[\group{z}]\ot_{\ql\clo}\qwc
 U)\\ =\Hom_{\ql\clo[H']}(V_\be\ot_{\ql\clo}V_\varpi,V_\be\ot_{\ql\clo}\qwc
 U)\,.\earr$$ Then $h\mapsto1\ot h$ is an isomorphism
 $\fV_\varpi\to\fV_{\be\ot\varpi}$ of vector spaces over $\qwc\Ga_k$ and
 one checks that $(1\ot h)x=1\ot h\cdot\be(x)$\,.

\sn The same argument, with $T,\Tr$ rather than $K_1,\Det$, yields
the commutativity of the right quadrangle, and the commutativity
of the bottom quadrangle follows from the formula for $\LL$ by
$\psi_l(\be)=\be^{\mr{Fr}}$. The diagram now implies that the top
quadrangle commutes.

\sn Recall that, [RW3, Proposition 11], the logarithm
$\Ll_\be:K_1(\lwb U)\to T(\qwb U)$ is integral for all $\be$. It
thus suffices to show that if $x\in T(\qwg)$ has image in $T\lwb
U)$ under the southwest arrow for every $\be$, then $x\in T(\lw
U)$.

\sn Letting $\be$ run through a set of representatives of the
$G(\ql\clo/\ql)$-action on the $\ql\clo$-irreducible characters of
$\group{z}$ we get an isomorphism
$\zl[\group{z}]\to\prod_\be\zl[\be]$. This induces isomorphisms
$T(\lwg)\to\bigoplus_\be T(\lwb U)$ and $T(\qwg)\to\bigoplus_\be
T(\qwb U)$. The first of these provides an $x'\in T(\lwg)$ with
the same images as $x$ for all $\be$, and the second gives
$x=x'\in T(\lwg)$\,.

\sn We now prove ii). The exact sequence defining $\fa$ gives the
top row since $\fa\subseteq\mr{rad}(\lwg)$\,, as $[G,G]$ is an
$l$-group. The bottom row is exact by ii) of Lemma {\sc e}. To see
the vertical surjectivity, write $\fu$ for the kernel of
$\defl_U^\uab: \lw U\onto\lw\uab$ (with $U=G[l]$); also write
$\fu^\be$ for the kernel of $\lwb U\onto\lwb\uab$ (with $\be$ as
before). Then the map $1+\fu^\be\to\tau_U(\fu^\be)$ induced by
$\Ll_\be:(\lwb U)\mal\to T(\lwb U)$ is surjective [RW3, 2b.~of
Proposition 11]. Identifying $\lwg$ and $\prod_\be\lwb U$ as in
the last paragraph of the proof of i) (also for the
abelianizations), and assembl-

\mn\bmp{8cm} ing our asserted diagram in terms of the
$\be$-decomposition, noting that the commutativity of the square
at right follows from that of the top quadrangle in i),
we deduce that $1+\fa\to\tau_G(\fa)$ is also surjective. \emp\hspace{23mm}\bmp{5cm}$\barr{ccccc}(\lwg)\mal&\sr{\be}{\to}&(\lwb U)\mal&&\\
\daz{\Ll}&&\daz{\Ll_\be}&&\\
T(\lwg)&\sr{\be}{\to}&T(\lwb U)&&\Box \earr$\emp

\bbn Recall that, for a pro-$l$ group $G=G(K/k)$ and an open subgroup
$G'\le G$,
$$\Res_G^\gs:\Hom^\ast(R_lG,\qwc\Ga_k)\to\Hom^\ast(R_l\gs,\qwc\Ga_{k'})$$
is defined in [RW8,\S1]. We partially extend this definition to
$l$-elementary $G$.

\Lemma{j} Let $G=\group{z}\times G[l]$ be $l$-elementary. If $G'$
is an open subgroup of $G$ containing $\group{z}$, define for
$f\in\Hom^\ast(R_lG,\qwc\Ga_k)$
$$\Res_G^\gs f=[\,\chi'\mapsto f(\ind_\gs^G\chi')+\sum_{r\ge1}\fra{\Psi^r}{l^r}(f(\psi_l^{r-1}\chi))\,]\in\Hom^\ast(R_lG',\qwc\Ga_{k'})\,,$$
  where $\chi'\in R_l\gs$\,, $\psi_l$ denotes the $l^{\mr{th}}$ Adams operation,
  $\chi\df\psi_l(\ind_\gs^G\chi')-\ind_\gs^G(\psi_l\chi')$ and $k'=K^\gs$. Then
$\Res$ is additive, integral and transitive. Moreover, the diagram
below, and so diagram (D2), commutes.
 {\small $$\barr{cccccccc}K_1(\lwg)&\sr{\Det}{\to}&\HOM(R_lG,(\lwc\Ga_k)\mal)&\sr{\LL}{\to}&\Hom^\ast(R_lG,\qwc\Ga_k)&\sr{\Tr}{\leftarrow}&T(\lwg)&\\
 \daz{\res_G^\gs}&&\daz{\res_G^\gs}&&\daz{\Res_G^\gs}&&\daz{\Res_G^\gs}&\\
 K_1(\lw\gs)&\sr{\Det}{\to}&\HOM(R_lG',(\lwc\Ga_{k'})\mal)&\sr{\LL}{\to}&\Hom^\ast(R_lG',\qwc\Ga_{k'})&\sr{\Tr}{\leftarrow}&T(\qw\gs)\earr$$}
\Stop
\proof We first observe that for a suitable power $l^{r_0}$
of $l$, $G^{l^{r_0}}\subset G'$, thus
$$\psi_l^{r_0-1}\chi=\psi_l^{r_0-1}(\psi_l(\ind_\gs^G\chi')-\ind_\gs^G(\psi_l\chi'))=0\,:$$ compare [loc.cit.]
and note, restricting attention to irreducible $\chi'$, that
$\chi'=\be'\ot\varpi'$ for some irreducible $\ql\clo$-characters $\be'$
and $\varpi'$ whose kernels contain $G'[l]\,,\,\group{z}$\,, respectively. Then
$\ind_\gs^G\chi'=\be\ot\ind_{G'}^{G}\varpi'$ with $\res_G^\gs\be=\be'$\,, whence
$$\psi_l(\ind_\gs^G\chi')-\ind_\gs^G(\psi_l\chi')=\psi_l(\be)\ot(\psi_l(\ind_{G'}^{G}\varpi')-\ind_{G'}^{G}(\psi_l\varpi'))\,.$$
For the first assertion, we just proceed as in [loc.cit.,
top of p.120 and Proof of Proposition A].

\sn Finally, the map $\Res_G^\gs$ on the very right is defined by
transporting
$\Res_G^\gs:\Hom^\ast(R_lG,\qwc\Ga_k)\to\Hom^\ast(R_l\gs,\qwc\Ga_{k'})$
to $T(\qwg)\to T(\qw\gs)$ by means of the isomorphism `$\Tr$'; in
particular, the right square commutes. The middle square commutes
because of the computations shown in [loc.cit., Proof of Lemma 1].
And commutativity of the left square is [RW2, Lemma 9]. Observing
that $\Ll=\Tr\me\LL\Det$\,, diagram (D2) is obtained from
arranging columns 1, 4 and 3 with $\Ll$ and $\Tr$ as horizontal
maps.  \qed

\bbn We close this section with adjusting the arguments in \S1 for
the proof of Proposition 3 for pro-$l$ groups to $l$-elementary
groups $G$.

\Lemma{k} \ben\item[i)] If $t_{K/k}\in T(\lwg)$ then there is a
unique
 torsion
$w \in \HOM(R_lG,(\lwc\Ga_k)\mal)$ with $wL_{K/k}\in\Det K_1(\lwg)$ and $\defl_G^\gab w=1$\,.
\item[ii)] Moreover, if $G$ has an abelian subgroup $\gs$ of index $l$, then $$w=1\iff\ver_\gab^\gs\la_{K^{[G,G]}/k}\equiv\la_{K/k'}\mod\tr_Q(\lw\gs)$$
where $k'=K^\gs$ and $Q=G/\gs$.
 \een\Stop
 \proof For i), using the diagram in ii) of Lemma {\sc i} to replace the one in the proof of [RW5, Proposition 2.2]
 there exists a $y\in(\lwg)\mal$ so $\Ll(y)=t_{K/k}$ and $\defl_G^\gab w=1$. Following the proof of
 [loc.cit., Proposition 2.4] one defines $w$ by $wL_{K/k}=\Det y$ and checks that $\defl_G^\gab w=1$ and $\LL(w)=0$.
 This implies $w$ is torsion by the indicated argument from [RW3, p.46] by observing that, while
 $\psi_l^n\chi$ is only a character $\be$ of $G/G[l]$ for large $n$, $\be=\infl_\gab^G\be'$ still implies
 that $w(\be)=(\defl_G^\gab w)(\be')=\chi(1)1$\,. The argument for the uniqueness of $w$ still works because
 [RW3, Lemma 12] is already proved for $\lwb(G[l])$ in the notation of Lemma {\sc i}.

 \sn More precisely, in the notation of the proof of ii) of Lemma {\sc i}, let $x\in 1+\fa$ have $\Det x$ torsion.
 By the commutativity of the left quadrangle in the proof of i) of Lemma {\sc i}, $\be(x)\in 1+\fu^\be$
 has $\Det \be(x)$ torsion, hence we have $\Det \be(x)=1$, and so it suffices to observe that
 $\HOM(R_lG,(\lwc\Ga_k)\mal)\to\prod_\be\HOM^\be(R_l U,(\lwc\Ga_k)\mal)$ is injective.

 \sn To verify ii), we follow the proof [RW5, equivalence of (1) and (2) in Proposition 3.2], except that
 we still need to show that the only torsion unit $e\in\lw\gs$ congruent to $1\mod\tr_Q(\lw\gs)$ is $1$, even when
 $G$ is $l$-elementary. Decomposing the torsion subgroup $H$ of $G$ as $H=\group{z}\times H'$ we have
 $G'=\Ga\times H'$ with $\Ga\simeq\Ga_{k'}$ and $\lw\gs=\prod_\be(\lwb\Ga)[H']$. Now $\be(e)\equiv 1\mod\tr_Q((\lwb\Ga)[H'])$,
 hence $\be(e)=1$ by Higman's theorem for $(\lwb\Ga)[H']$, see [RW3, p.47]. This holds for all $\be$, hence $e=1$. \qed

 \def\kq{{\ol K}}
 \def\km{{\ol M}}
 \def\scd{{\mr{scd}}}
 \def\cor{{\mr{cor}}}
 \def\cd{{\mr{cd}}}

 \bbn{\sc Appendix}

 \bn The proof of [RW2, Proposition 12 (a)] refers to [RWt,
 Proposition 4.8] which, however, requires Leopoldt's conjecture
 (see [RWt, Lemma 3.4]). We recall the statement made in [RW2]
 (suppressing the index $\infty$ on $K$ and $G$ as well
 as the $\ti{\pht{\mho}}$ on $\mho_S$)\,:
 \bqo {\em If $N$ is a finite normal subgroup of $G$ with factor group
 $\gq$ and fixed field $\kq=K^N$, then \ $\defl_G^{\gq}:
 K_0T(\La G)\to K_0T(\La\gq)$ \ takes $\mho_S=\mho_{K/k,S}$ to
 $\ol\mho_S=\mho_{\kq/k,S}$.}\eqo The refinement $\mho_{K/k,S}$ of the Iwasawa module $X=X_{K/k,S}\df G(M/K)$,
 with $M$ the maximal abelian $S$-ramified $l$-extension of $K$, is described in
 [RW1,\S1].

 \mn Here is a direct argument for the above claim.

 \sn  Let $\ol M$ be the maximal abelian $S$-ramified
 $l$-extension of $\kq$, hence $G(\ol M/\kq)$ is the biggest abelian pro-$l$ quotient of $G(M/\kq)$. Consider the diagram below,
  where $\ti X$ is the pushout of $i$ and $\ver$, and the transfer $G(M/\kq)\to X^N$ factors through $G(\ol M/\kq)$ since $X^N$ is abelian pro-$l$.
 The bottom row is called the deflation of the top one in [RWt,\S3.2].


 $$\barr{ccccc}X=G(M/K)&\into&G(M/k)&\onto&G(K/k)=G\\
 \cap&&\pl&&\sda\\
 G(M/\ol K)&\into&G(M/k)&\onto&G(\ol K/k)=\ol G\\
 \sda&&\sda&&\pl\\
 \ol X=G(\ol M/\ol K)&\sr{i}{\into}&G(\ol M/\kq)&\onto&\ol G\\
 \daz{\ver}&&\da&&\pl\\
 X^N&\into&\ti X&\onto&\ol G\,,\earr$$
  By the Appendix 4.A analogue of Lemma 3.2 [loc.cit.] the translation functor turns the bottom two
  rows into \footnote{here $\De G$ denotes the augmentation ideal of $\La G$}
 \nf{A1}{\barr{ccccc}\ol X&\into&\ol Y&\onto&\De\gq\\
 \da&&\da&&\pl\\ X^N&\into&Y_N&\onto&(\De G)_N\,,\earr}
 by replacing the bottom one by the equivalent extension given by the $\La G$-analogue of Lemma 3.3.
 Here we should note that $Y$ has projective dimension $\le1$ over $\La G$ [RW1, Theorem 1],
 hence $Y$ has projective dimension $\le1$ over $\zl[N]$ [loc.cit., proof of 2.~of Proposition 4, which does not need $M$ to be finite].

 \sn Suppose that we know that $\ver:G(\ol M/\kq)\to X^N$ is an isomorphism,
 hence that the extensions in (A1) are equivalent. If we use,  [loc.cit.,\S1],
 $X\into Y\to\La G\onto\zl$ to compute $\mho$ then taking $N$-coinvariants computes
 the analogous $\defl_G^\gq(\mho)$ for the bottom row of (A1) (cf.~the analogy
 with [RWt, Lemma 4B.1, p.41]). By (A1), the same procedure for the top row computes $\ol\mho$. Thus $\defl_G^\gq(\mho)=\ol\mho$.

 \sn Concerning $G(\ol M/\ol K)\sr{\ver}{\lto}X^N$\,, let $L$ be the maximal $S$-ramified Galois extension of $k$
 or, equivalently, of $k_\infty$. Denote the corresponding Galois
 groups by $\fG$ and $\fH$, respectively; so
 $\fH\into\fG\onto\Ga_k$ is exact. Moreover, set $V=G(L/K)$ and $U=G(L/\ol
 K)$, whence $V\into U\onto N$.

 \sn Assume that we already know
 $\scd_l(\fH)=2$. Then we proceed as follows. As $U$ is open in
 $\fH$, it follows that also $\scd_l(U)=2$ by [NSW, p.139/140]. The proof of
 [\,(i) $\implies$ (ii)\,] of [loc.cit., Theorem 3.6.4, p.160] gives the isomorphism \
 $H^2(V,\Z)(l)_N\sr{\cor}{\lto}H^2(U,\Z)(l)$ (see [loc.cit., 3.3.8, p.142])
 and so
 $(U^\ab)_l\sr{\ver}{\lto}(V^\ab)_l^N$. Since
 $(U^\ab)_l=G(\ol M/\ol K)=\ol X$\,,  $(V^\ab)_l=G(M/K)=X$ finishes the proof.

 \sn Hence it remains to show $\scd_l(\fH)=2$. Now, $\scd_l(\fH)\le2$ is a consequence of the weak Leopoldt
 conjecture (see [loc.cit., 10.3.26, p.549]) and then $\scd_l(\fH)=2$ results from the
 remark following it, of which we add a proof\,:
  \ Assume $\scd_l(\fH)\le1$. Then
 $\cd_l(\fH)\le1$ and $\cd_l(\Ga_k)=1$\,, hence
 $2=\cd_l(\fG)\le\cd_l(\Ga_k)+\cd_l(\fH)$ implies $\cd_l(\fH)=1$.
 Note that $\cd_l(\fG)=2$ by [loc.cit., 10.9.3, p.587].
  Denoting a Sylow-$l$ subgroup of $\fH$ by
 $\fH_l$, we have $\cd_l(\fH_l)=1=\scd_l(\fH_l)$ [loc.cit.,3.3.6, p.141]
 and thus $H^2(\fH_l,\Z)(l)=0$ [loc.cit.,3.3.4, p.139]. Hence
 $H^1(\fH_l,\ql/\zl)=0$ and so $\fH_l=1$, contradicting
 $\cd_l(\fH_l)=1$.

 \bbbn{\sc References}
 \small
 \bbn\btb{rp{13cm}}
 \,[CR] &  Curtis, C.W.~and Reiner, I., {\em Methods of Representation Theory,
           vol. 2.} John Wiley \& Sons (1987) \\
 \,[DR]&   Deligne, P.~and Ribet, K., {\em Values of abelian
           $L$-functions at negative integers
           over totally real fields.} Invent.~Math. {\bf 59} (1980), 227-286\\
 \,[FK]   & Fukaya, T., Kato, K., {\em A formulation of conjectures on $p$-adic zeta
            functions in noncommutative Iwasawa theory.} Proceedings of
            the St.~Petersburg Mathematical Society, vol.~XII (ed.
            N.N.~Uraltseva), AMS Translations -- Series 2, {\bf 219} (2006),
            1-86\\
 \,[HIÖ] & Hawkes, T., Isaacs, I.M.~and Özaydin, M.\,, {\em On the Möbius function of a finite group.} Rocky Mountain J.~of
           Mathematics {\bf 19} (1989), 1003-1033 \\
\etb
\noi \btb{rp{13cm}}
 \,[Ka]  &  Kato, K.\,, {\em Iwasawa theory and generalizations.}
           Proc.\,ICM, Madrid, Spain, 2006; European
           Math.\,Soc.\,(2007), 335-357\\
 \,[La]  & Lau, I., {\em Algebraic contributions to equivariant Iwasawa theory.} Ph.D.\,thesis, Universität Augsburg (2010)\\
 \,[NSW] & Neukirch, J., Schmidt, A.~and Wingberg, K., {\em Cohomology of
           Number Fields.} Springer Grundlehren der math.~Wiss.~{\bf 323} (2000)\\
 \,[RW] &  Ritter, J.~and Weiss, A.\,, \newline
           \hsp{2}t. {\em The Lifted Root Number Conjecture and Iwasawa theory.}
           Memoirs of the AMS \hsp{6}{\bf 157/748} (2002)\newline
           \hsp{2}1. {\em Towards equivariant Iwasawa theory.} manuscripta math.~{\bf 109} (2002), 131-146\newline
           \hsp{2}2. {\em $\cdots$, II.}
           Indag.\,Mathemat.~{\bf 15} (2004), 549-572\newline
           \hsp{2}3. {\em $\cdots$, III.} Math.\,Ann.~{\bf 336}
           (2006), 27-49\newline
           \hsp{2}4. {\em $\cdots$, IV.} Homology, Homotopy and
           Applications {\bf 7} (2005), 155-171\newline
           \hsp{2}5. {\em Non-abelian pseudomeasures and congruences
           between abelian Iwasawa $L$- \hsp{6}functions.}
           Pure and Applied Math.~Quarterly {\bf 4} (2008), 1085-1106\newline
           \hsp{2}6. {\em The integral logarithm in Iwasawa
           theory: an exercise.} To appear in Journal de \hsp{6}Th\'eorie des Nombres Bordeaux {\bf 22} (2010)\newline
           \hsp{2}7. {\em Congruences between abelian pseudomeasures.}
           Math.Res.\,Lett.~{\bf 15} (2008), 715-725\newline
           \hsp{2}8. {\em Equivariant Iwasawa Theory: An Example.}
           Documenta Math.~{\bf 13} (2008), 117-129\newline
           \hsp{2}9. {\em Congruences between abelian pseudomeasures, II.} arXiv:1001.2091
           [math.NT]\\
 \,[Se]&   Serre, J.-P.\,, {\em Sur le r\'esidu de la fonction
           z\^{e}ta $p$-adique d'un corps de nombres.}\newline
           C.R.Acad.Sci.~Paris {\bf 287} (1978), s\'erie A,
           183-188\\
 \,[Wi]   & Wiles, A., {\em The Iwasawa conjecture for totally real
            fields.} Annals of Math. {\bf 131} (1990), 493-540
 \etb
 \tiny
 \vspace*{1.5cm}
 \bct Jürgen Ritter, Schnurbeinstraße 14, 86391 Deuringen,
 Germany; {\tt jr@ritter-maths.de}\\
 Alfred Weiss, Department of Mathematics, University of Alberta,
 Edmonton, AB, Canada T6G 2G1; {\tt weissa@ualberta.ca}\ect

 \end{document}